\documentclass[12pt,a4paper]{amsart}
\usepackage{amsfonts}
\usepackage{amsmath}
\usepackage{amssymb}

\def\gH{\mathfrak{H}}

\def\gM{\mathfrak{M}}
\def\gN{\mathfrak{N}}

\begin{document}

\title[On Lesieur's Measured Quantum Groupoids] 
{On Lesieur's Measured Quantum Groupoids}
\author{Michel Enock}
\address{Institut de Math\'ematiques de Jussieu, Unit\'{e} Mixte Paris 6 / Paris 7 /
CNRS de Recherche 7586 \\175, rue du Chevaleret, Plateau 7E, F-75013 Paris}
 \email{enock@math.jussieu.fr}
\date{june 07}

\begin{abstract}
In his thesis ([L1]), which is published in an expended and revised version ([L2]), Franck Lesieur had introduced a notion of measured quantum groupoid, in the setting of von Neumann algebras, using intensively the notion of pseudo-multiplicative unitary, which had been introduced in a previous article of the author, in collaboration with Jean-Michel Vallin [EV]. In [L2], the axioms given are very complicated and are here simplified.  \end{abstract}

\maketitle
\newpage
\section{Introduction}
\label{intro}
\subsection{}
 In two articles ([Val1], [Val2]), J.-M. Vallin has introduced two notions (pseudo-multiplicative
unitary, Hopf-bimodule), in order to generalize, up to the groupoid
case, the classical notions of multiplicative unitary [BS] and of Hopf-von Neumann algebras [ES]
which were introduced to describe and explain duality of groups, and leaded to appropriate notions
of quantum groups ([ES], [W1], [W2], [BS], [MN], [W3], [KV1], [KV2], [MNW]). 
\\ In another article [EV], J.-M. Vallin and the author have constructed, from a depth 2 inclusion of
von Neumann algebras $M_0\subset M_1$, with an operator-valued weight $T_1$ verifying a regularity
condition, a pseudo-multiplicative unitary, which leaded to two structures of Hopf bimodules, dual
to each other. Moreover, we have then
constructed an action of one of these structures on the algebra $M_1$ such that $M_0$
is the fixed point subalgebra, the algebra $M_2$ given by the basic construction being then
isomorphic to the crossed-product. We construct on $M_2$ an action of the other structure, which
can be considered as the dual action.
\\  If the inclusion
$M_0\subset M_1$ is irreducible, we recovered quantum groups, as proved and studied in former papers
([EN], [E1]).
\\ Therefore, this construction leads to a notion of "quantum groupo\"{\i}d", and a construction of a
duality within "quantum groupo\"{\i}ds". 
\subsection{}
In a finite-dimensional setting, this construction can be
mostly simplified, and is studied in [NV1], [BSz1],
[BSz2], [Sz], [Val3], [Val4], and examples are described. In [NV2], the link between these "finite quantum
groupo\"{\i}ds" and depth 2 inclusions of
$II_1$ factors is given. 
\subsection{}
Franck Lesieur introduced in his thesis [L1] a notion of "measured quantum groupoids", in which a modular hypothesis on the basis is required. Mimicking in a wider setting the technics of Kustermans and Vaes [KV], he obtained then a pseudo-multiplicative unitary, which, as in the quantum group case, "contains" all the information of the object (the von Neuman algebra, the coproduct, the antipod, the co-inverse). Unfortunately, the axioms chosen then by Lesieur don't fit perfectely with the duality (namely, the dual object doesnot fit the modular condition on the basis chosen in [L1]), and, for this purpose, Lesieur gave the name of "measured quantum groupoids" to a wider class [L2], whose axioms could be described as the analog of [MNW], in which a duality is defined and studied, the initial objects considered in [L1] being denoted now "adapted measured quantum groupoids". In [E3] had been shown that, with suitable conditions, the objects constructed in [EV] from depth 2 inclusions, are "measured quantum groupoids" in this new setting. 
\subsection{}
Unfortunately, the axioms given in ([L2], 4) are very complicated, and there was a serious need for simplification. This is the goal of this article. 
 \subsection{}
This article is organized as follows : 
\newline
In chapter \ref{pr} are recalled all the definitions and constructions needed for that theory, namely Connes-Sauvageot's relative tensor product of Hilbert spaces, fiber product of von Neumann algebras, and Vaes' Radon-Nikodym theorem. 
\newline
The chapter \ref{pmu} is a r\'esum\'e of Lesieur's basic result ([L2], 3), namely the construction of a pseudo-multiplictaive unitary associated to a Hopf-bimodule, when exist a left-invariant operator-valued weight, and a right-invariant valued weight. 
\newline
The chapter \ref{coinverse} is mostly inspired from Lesieur's "adapted measured quantum groupoids" ([L2], 9), with a wider hypothesis, namely, that there exists a weight on the basis such that the modular automorphism groups of two lifted weights (via the two operator-valued weights) commute. This hypothesis allows us to use Vaes' theorem, and is a nice generalization of the existence of a relatively invariant measure on the basis of a groupoid. With that hypothesis, mimicking ([L2], 9), we construct a co-inverse and a scaling group. 
\newline
In chapter \ref{AGB}, we go on with the same hypothesis. It allows us to construct two automorphism groups on the basis, which appear to be invariant under the relatively invariant weight introduced in chapter \ref{coinverse}
\newline
It is then straightforward to get that we are now in pr\'esence of Lesieur's "measured quantum groupoids" ([L2], 4) and chapter \ref{MQG} is devoted to main properties of these. 
 \subsection{}
 The author is indebted to Frank Lesieur, Stefaan Vaes, Leonid Va\u{\i}nerman, and especially Jean-Michel Vallin, for many fruitful conversations. 

\section{Preliminaries}
\label{pr}
In this chapter are mainly recalled definitions and notations about Connes' spatial
theory (\ref{spatial}, \ref{rel}) and the fiber product construction (\ref{fiber}, \ref{slice})
which are the main technical tools of the theory of measured quantum theory.

\subsection{Spatial theory [C1], [S2], [T]}
 \label{spatial}
 Let $N$ be a von Neumann algebra, 
and let $\psi$ be a faithful semi-finite normal weight on $N$; let $\gN_{\psi}$, 
$\gM_{\psi}$, $H_{\psi}$, $\pi_{\psi}$, $\Lambda_{\psi}$,$J_{\psi}$, 
$\Delta_{\psi}$,... be the canonical objects of the Tomita-Takesaki construction 
associated to the weight $\psi$. Let $\alpha$  be a non-degenerate normal representation of $N$ on a
Hilbert space
$\mathcal{H}$. We may as well consider $\mathcal{H}$ as a left $N$-module, and write it then
$_\alpha\mathcal{H}$. Following ([C1], definition 1), we define the set of 
$\psi$-bounded elements of $_\alpha\mathcal{H}$ as :
\[D(_\alpha\mathcal{H}, \psi)= \lbrace \xi \in \mathcal{H};\exists C < \infty ,\| \alpha (y) \xi\|
\leq C \| \Lambda_{\psi}(y)\|,\forall y\in \gN_{\psi}\rbrace\]
Then, for any $\xi$ in $D(_\alpha\mathcal{H}, \psi)$, there exists a bounded operator
$R^{\alpha,\psi}(\xi)$ from $H_\psi$ to $\mathcal{H}$,  defined, for all $y$ in $\gN_\psi$ by :
\[R^{\alpha,\psi}(\xi)\Lambda_\psi (y) = \alpha (y)\xi\]
This operator belongs to $Hom_N (H_\psi , \mathcal{H})$;
therefore, for any
$\xi$, $\eta$ in
$D(_\alpha\mathcal{H}, \psi)$, the operator :
\[\theta^{\alpha,\psi} (\xi ,\eta ) =  R^{\alpha,\psi}(\xi)R^{\alpha,\psi}(\eta)^*\]
belongs to $\alpha (N)'$; moreover, $D(_\alpha\mathcal{H}, \psi)$ is dense ([C1], lemma 2), stable
under $\alpha (N)'$, and the linear span generated by the operators $\theta^{\alpha,\psi} (\xi ,\eta
)$ is a weakly dense ideal in $\alpha (N)'$. 
 \newline
 With the same hypothesis, the operator :
\[<\xi,\eta>_{\alpha,\psi} = R^{\alpha,\psi}(\eta)^* R^{\alpha,\psi}(\xi)\]
belongs to $\pi_{\psi}(N)'$. Using Tomita-Takesaki's theory, this last algebra is equal to
$J_\psi
\pi_{\psi}(N)J_\psi$, and therefore anti-isomorphic to $N$ (or isomorphic to the opposite von
Neumann algebra $N^o$). We shall consider now $<\xi,\eta>_{\alpha,\psi}$ as an element of $N^o$, and
the linear span generated by these operators is a dense algebra in $N^o$. More precisely ([C], lemma 4, and [S2], lemme 1.5b), we get that $<\xi, \eta>_{\alpha, \psi}^o$ belongs to $\gM_\psi$, and that :
 \[\Lambda_{\psi}(<\xi, \eta>_{\alpha, \psi}^o)=J_\psi R^{\alpha, \psi}(\xi)^*\eta\]
 \newline
 If $y$ in $N$ is analytical with respect to $\psi$, and if $\xi\in D(_\alpha\mathcal{H}, \psi)$, then we get that $\alpha(y)\xi$ belongs to $D(_\alpha\mathcal{H}, \psi)$ and that :
 \[R^{\alpha,\psi}(\alpha(y)\xi)=R^{\alpha,\psi}(\xi)J_\psi\sigma_{-i/2}^{\psi}(y^*)J_\psi\]
  So, if $\eta$ is another $\psi$-bounded element of $_\alpha\mathcal H$, we get :
 \[<\alpha(y)\xi, \eta>_{\alpha, \psi}^o=\sigma_{i/2}^\psi(y)<\xi, \eta>_{\alpha, \psi}^o\]
There exists ([C], prop.3) a family  $(e_i)_{i\in I}$ of 
$\psi$-bounded elements of $_\alpha\mathcal H$, such that
\[\sum_i\theta^{\alpha, \psi} (e_i ,e_i )=1\]
Such a family will be called an $(\alpha,\psi )$-basis of $\mathcal H$. 
 \newline
 It is possible ([EN] 2.2) to construct 
an $(\alpha,\psi )$-basis of $\mathcal H$, $(e_i)_{i\in I}$, such that the
operators $R^{\alpha, \psi}(e_i)$ are partial isometries with final supports 
$\theta^{\alpha, \psi}(e_i ,e_i )$ 2 by 2 orthogonal, and such that, if $i\neq j$, then 
$<e_i ,e_j>_{\alpha, \psi}=0$. Such a family will be called an $(\alpha, \psi)$-orthogonal basis of $\mathcal H$. 
\newline
 We have, then :
 \[R^{\alpha, \psi}(\xi)=\sum_i\theta^{\alpha, \psi}(e_i, e_i)R^{\alpha, \psi}(\xi)=\sum_iR^{\alpha, \psi}(e_i)<\xi, e_i>_{\alpha, \psi}\]
 \[<\xi, \eta>_{\alpha, \psi}=\sum_i<\eta, e_i>_{\alpha, \psi}^*<\xi, e_i>_{\alpha, \psi}\]
 \[\xi=\sum_iR^{\alpha, \psi}(e_i)J_\psi\Lambda_\psi(<\xi, e_i>^o_{\alpha, \psi})\]
 the sums being weakly convergent. Moreover, we get that, for all $n$ in $N$, $\theta^{\alpha, \psi}(e_i, e_i)\alpha(n)e_i=\alpha(n)e_i$, and $\theta^{\alpha, \psi}(e_i, e_i)$ is the orthogonal projection on the closure of the subspace $\{\alpha(n)e_i, n\in N\}$. 
 \newline
 If $\theta\in Aut N$, then it is straightforward to get that $D(_{\alpha\circ\theta}\mathcal H, \psi\circ\theta)=D(_\alpha\mathcal H, \psi)$, and then, we get that, for any $\xi$, $\eta$ in $D(_\alpha\mathcal H, \psi)$ :
 \[<\xi, \eta>_{\alpha\circ\theta, \psi\circ\theta}^o=\theta^{-1}(<\xi, \eta>_{\alpha, \psi}^o)\]
Let $\beta$ be a normal non-degenerate
anti-representation
 of
$N$ on
$\mathcal{H}$. We may then as well consider $\mathcal{H}$ as a right $N$-module, and write it $\mathcal{H}_\beta$, or
consider
$\beta$ as a normal non-degenerate representation of the opposite von Neumann algebra $N^o$, and
consider
$\mathcal{H}$ as a left $N^o$-module. 
\newline
We can then define on $N^o$ the opposite faithful
 semi-finite normal weight $\psi ^o$; we have $\gN_{\psi ^o}=\gN_\psi^*$, and 
 the Hilbert space $H_{\psi ^o}$ will be, as usual, identified with $H_\psi$, 
by the identification, for all $x$ in $\gN_\psi$, of 
$\Lambda_{\psi ^o}(x^*)$ with $J_\psi \Lambda_\psi (x)$.
\newline
From these remarks, we infer that the set of 
$\psi ^o$-bounded elements of
$\mathcal{H}_\beta$ is :
\[D(\mathcal{H}_\beta, \psi^o) = \lbrace\xi\in\mathcal{H} ;\exists C < \infty ,
\|\beta (y^*)\xi\|
\le C \| \Lambda_{\psi}(y)\|,\forall y\in \gN_{\psi}\rbrace\]
and, for any $\xi$ in $D(\mathcal{H}_\beta, \psi^o)$ and $y$ in $\gN_\psi$,
 the bounded operator $R^{\beta,\psi^{o}}(\xi )$ is given by the formula :
\[R^{\beta,\psi{^o}}(\xi)J_{\psi}\Lambda_\psi (y) = \beta (y^*)\xi\]
This operator belongs to $Hom_{N^{o}}(H_\psi ,\mathcal{H} )$.
Moreover, $D(\mathcal{H}_\beta, \psi^o)$ is dense, stable under $\beta( N)'=P$, 
and, for all $y$ in $P$, we have :
\[R^{\beta,\psi{^o}}(y\xi)= yR^{\beta,\psi{^o}}(\xi)\]
Then, for any $\xi$, $\eta$ in $D(\mathcal{H}_\beta, \psi^o)$, the operator 
\[\theta^{\beta, \psi^{o}}(\xi ,\eta )=R^{\beta,\psi{^o}}(\xi)R^{\beta,\psi{^o}}(\eta)^*\] 
belongs to $P$, and the linear span generated by these operators is a dense ideal in $P$;
moreover, the operator-valued product 
$<\xi ,\eta >_{\beta,\psi^o}= R^{\beta,\psi{^o}}(\eta)^*R^{\beta,\psi{^o}}(\xi)$
belongs to $\pi_\psi (N)$; we shall consider now, for simplification, that $<\xi ,\eta
>_{\beta,\psi^o}$ belongs to $N$, and the linear span generated by these operators is a dense algebra
in $N$, stable under multiplication by analytic elements with respect to $\psi$. More precisely, $<\xi ,\eta
>_{\beta,\psi^o}$ belongs to $\gM_\psi$ ([C], lemma 4) and we have ([S1], lemme 1.5) 
\[\Lambda_\psi(<\xi, \eta>_{\beta,\psi^o})=R^{\beta,\psi{^o}}(\eta)^*\xi\]
A $(\beta,\psi^o)$-basis of $\mathcal H$ is a family 
$(e_i )_{i\in I}$ of $\psi^o$-bounded elements of $\mathcal H_\beta$, such that 
\[\sum_i\theta^{\beta, \psi^{o}}(e_i ,e_i )=1\]
We have then, for all $\xi$ in $D(\mathcal H_\beta)$ :
\[\xi=\sum_i R^{\beta,\psi^o}(e_i)\Lambda_\psi(<\xi, e_i>_{\beta,\psi^o})\]
It is possible to choose the $(e_i )_{i\in I}$ such that
the $R^{\beta, \psi{^o}}(e_i)$ are partial isometries, with final supports
$\theta^{\beta, \psi^{o}}(e_i ,e_i )$ 2 by 2 orthogonal, and such that $<e_i, e_j>_{\beta, \psi^o}=0$ if $i\neq j$; such
a family will be then called a $(\beta, \psi^o)$-orthogonal basis of $\mathcal H$. We have then 
\[R^{\beta, \psi{^o}}(e_i)=\theta^{\beta, \psi^{o}}(e_i ,e_i )R^{\beta, \psi{^o}}(e_i)=
R^{\beta, \psi{^o}}(e_i)<e_i, e_i>_{\beta, \psi^o}\]
 Moreover, we get that, for all $n$ in $N$, and for all $i$, we have :
\[\theta^{\beta, \psi^{o}}(e_i ,e_i )\beta (n)e_i=\beta (n)e_i\]
and that  $\theta^{\beta, \psi^{o}}(e_i ,e_i )$ is the orthogonal projection on the closure of the subspace $\{\beta (n)e_i, n\in N\}$.

\subsection{Jones' basic construction and operator-valued weights}
\label{basic}
Let $M_0\subset M_1$ be an inclusion of von Neumann algebras  (for simplification, these algebras will be supposed to be $\sigma$-finite), equipped with a normal faithful semi-finite operator-valued weight $T_1$ from $M_1$ to $M_0$ (to be more precise, from $M_1^{+}$ to the extended positive elements of $M_0$ (cf. [T] IX.4.12)). Let $\psi_0$ be a normal faithful semi-finite weight on $M_0$, and $\psi_1=\psi_0\circ T_1$; for $i=0,1$, let $H_i=H_{\psi_i}$, $J_i=J_{\psi_i}$, $\Delta_i=\Delta_{\psi_i}$ be the usual objects constructed by the Tomita-Takesaki theory associated to these weights. Following ([J], 3.1.5(i)), the von Neumann algebra $M_2=J_1M'_0J_1$ defined on the Hilbert space $H_1$ will be called the basic construction made from the inclusion $M_0\subset M_1$. We have $M_1\subset M_2$, and we shall say that the inclusion $M_0\subset M_1\subset M_2$ is standard.  
 \newline
Following ([EN] 10.6), for $x$ in $\gN_{T_1}$, we shall define $\Lambda_{T_1}(x)$ by the following formula, for all $z$ in $\gN_{\psi_{0}}$ :
\[\Lambda_{T_1}(x)\Lambda_{\psi_{0}}(z)=\Lambda_{\psi_1}(xz)\]
Then, $\Lambda_{T_1}(x)$ belongs to $Hom_{M_{0}^o}(H_{0}, H_1)$; if $x$, $y$ belong to $\gN_{T_1}$, then $\Lambda_{T_1}(x)^*\Lambda_{T_1}(y)=T_1(x^*y)$, and $\Lambda_{T_1}(x)\Lambda_{T_1}(y)^*$ belongs to $M_{2}$.
\newline
 Using then Haagerup's construction ([T], IX.4.24), it is possible to construct a normal semi-finite faithful operator-valued weight $T_2$ from $M_2$ to $M_1$ ([EN], 10.7), which will be called the basic construction made from $T_1$. If $x$, $y$ belong to $\gN_{T_1}$, then $\Lambda_{T_1}(x)\Lambda_{T_1}(y)^*$
belongs to $\gM_{T_{2}}$, and $T_{2}(\Lambda_{T_1}(x)\Lambda_{T_1}(y)^*)=xy^*$. 
\newline
By Tomita-Takesaki theory, the Hilbert space $H_1$ bears a natural structure of $M_1-M_1^o$-bimodule, and, therefore, by restriction, of $M_0-M_0^o$-bimodule. Let us write $r$ for the canonical representation of $M_0$ on $H_1$, and $s$ for the canonical antirepresentation given, for all $x$ in $M_0$, by $s(x)=J_1r(x)^*J_1$. Let us have now a closer look to the subspaces $D(H_{1s}, \psi_0^o)$ and $D(_rH_1, \psi_0)$. If $x$ belongs to $\gN_{T_1}\cap\gN_{\psi_1}$, we easily get that $J_1\Lambda_{\psi_1}(x)$ belongs to $D(_rH_1, \psi_0)$, with :
\[R^{r, \psi_0}(J_1\Lambda_{\psi_1}(x))=J_1\Lambda_{T_1}(x)J_0\]
and $\Lambda_{\psi_1}(x)$ belongs to $D(H_{1s}, \psi_0)$, with :
\[R^{s, \psi_0^o}(\Lambda_{\psi_1}(x))=\Lambda_{T_1}(x)\]
In ([E3], 2.3) was proved that the subspace $D(H_{1s}, \psi_0^o)\cap D(_rH_1, \psi_0)$ is dense in $H_1$; let us write down and precise this result :

\subsubsection{{\bf Proposition}}
\label{propbasic}
{\it  Let us keep on the notations of this paragraph; let $\mathcal T_{\psi_1, T_1}$ be the algebra made of elements $x$ in $\gN_{\psi_1}\cap\gN_{T_1}\cap\gN_{\psi_1}^*\cap\gN_{T_1}^*$, analytical with respect to $\psi_1$, and such that, for all $z$ in $\mathbb{C}$, $\sigma^{\psi_1}_z(x)$ belongs to $\gN_{\psi_1}\cap\gN_{T_1}\cap\gN_{\psi_1}^*\cap\gN_{T_1}^*$. Then :
\newline
(i) the algebra $\mathcal T_{\psi_1, T_1}$ is weakly dense in $M_1$; it will be called Tomita's algebra with respect to $\psi_1$ and $T_1$; 
\newline
(ii) for any $x$ in  $\mathcal T_{\psi_1, T_1}$, $\Lambda_{\psi_1}(x)$ belongs to $D(H_{1s}, \psi_0^o)\cap D(_rH_1, \psi_0)$;
\newline
(iii) for any $\xi$ in $D(H_{1s}, \psi_0^o))$, there exists a sequence $x_n$ in $\mathcal T_{\psi_1, T_1}$ such that $\Lambda_{T_1}(x_n)=R^{s, \psi_0^o}(\Lambda_{\psi_1}(x))$ is weakly converging to $R^{s, \psi_0^o}(\xi)$ and $\Lambda_{\psi_1}(x_n)$ is converging to $\xi$. }

\begin{proof}
The result (i) is taken from ([EN], 10.12); we get in ([E3], 2.3) an increasing sequence of projections $p_n$ in $M_1$, converging to $1$, and elements $x_n$ in $\mathcal T_{\psi_1, T_1}$ such that $\Lambda_{\psi_1}(x_n)=p_n\xi$. So, (i) and (ii) were obtained in ([E3], 2.3) from this construction. More precisely, we get that :
\begin{eqnarray*}
T_1(x_n^*x_n)&=&<R^{s, \psi^o}(\Lambda_{\psi_1}(x_n)), R^{s, \psi_0^o}(\Lambda_{\psi_1}(x_n))>_{s, \psi_0^o}\\
&=&<p_n\xi, p_n\xi>_{s, \psi_0^o}\\
&=&R^{s, \psi^o}(\xi)^*p_nR^{s, \psi^o}(\xi)
\end{eqnarray*}
which is increasing and weakly converging to $<\xi, \xi>_{s, \psi_0^o}$.  \end{proof}

We finish by writing a proof of this useful lemma, we were not able to find in litterature   :
\subsubsection{{\bf Lemma}}
\label{lemT}
{\it Let $M_0\subset M_1$ be an inclusion of von neumann algebras, equipped with a normal faithful semi-finite operator-valued weight $T$ from $M_1$ to $M_0$. Let $\psi_0$ be a normal semi-finite faithful weight on $M_0$, and $\psi_1=\psi_0\circ T$; if $x$ is in $\gN_T$, and if $y$ is in $M'_0\cap M_1$, analytical with respect to $\psi_1$, then $xy$ belongs to $\gN_T$. }
\begin{proof}
Let $a$ be in $\gN_{\psi_0}$; then $xa$ belongs to $\gN_{\psi_1}$, and $xya=xay$ belongs to $\gN_{\psi_1}$; moreover, let us consider the element $T(y^*x^*xy)$ of the positive extended part of $M_0^+$; we have :
\begin{eqnarray*}
<T(y^*x^*xy), \omega_{\Lambda_{\psi_0}}(a)>=\psi_1(a^*y^*x^*xya)=\|\Lambda_{\psi_1}(xay)\|^2=\\
=\|J_{\psi_1}\sigma_{-i/2}^{\psi_1}(y^*)J_{\psi_1}\Lambda_{\psi_1}(xa)\|^2=\|J_{\psi_1}\sigma_{-i/2}^{\psi_1}(y^*)J_{\psi_1}\Lambda_T(x)\Lambda_{\psi_0}(a)\|^2
\end{eqnarray*}
from which we get that $T(y^*x^*xy)$ is bounded and \[T(y^*x^*xy)\leq \|\sigma_{-i/2}^{\psi_1}(y^*)\|^2T(x^*x)\] \end{proof}


\subsection{Relative tensor product [C1], [S2], [T]}
\label{rel}
Using the notations of \ref{spatial}, let now $\mathcal{K}$ be another Hilbert space on which there exists
a non-degenerate representation
$\gamma$ of
$N$. Following J.-L. Sauvageot ([S2], 2.1), we define
the relative tensor product $\mathcal{H}\underset{\psi}{_\beta\otimes_\gamma}\mathcal{K}$ as the
Hilbert space obtained from the algebraic tensor product $D(\mathcal{H}_\beta ,\psi^o )\odot
\mathcal{K} $ equipped with the scalar product defined, for $\xi_1$, $\xi_2$ in $D(\mathcal{H}_\beta
,\psi^o )$,
$\eta_1$, $\eta_2$ in $\mathcal{K}$, by 
\[(\xi_1\odot\eta_1 |\xi_2\odot\eta_2 )=(\gamma(<\xi_1 ,\xi_2 >_{\beta,\psi^o})\eta_1 |\eta_2 )\]
where we have identified $N$ with $\pi_\psi (N)$ to simplifly the notations.
\newline
The image of $\xi\odot\eta$ in $\mathcal{H}\underset{\psi}{_\beta\otimes_\gamma}\mathcal{K}$ will be
denoted by
$\xi\underset{\psi}{_\beta\otimes_\gamma}\eta$. We shall use intensively this construction; one
should bear in mind that, if we start from another faithful semi-finite normal weight $\psi '$, we
get another Hilbert space $\mathcal{H}\underset{\psi'}{_\beta\otimes_\gamma}\mathcal{K}$; there exists an isomorphism $U^{\psi, \psi'}_{\beta, \gamma}$ from $\mathcal{H}\underset{\psi}{_\beta\otimes_\gamma}\mathcal{K}$ to $\mathcal{H}\underset{\psi'}{_\beta\otimes_\gamma}\mathcal{K}$, which is unique up to some functorial property ([S2], 2.6) (but this isomorphism
does not send  $\xi\underset{\psi}{_\beta\otimes_\gamma}\eta$ on
$\xi\underset{\psi'}{_\beta\otimes_\gamma}\eta$ !). 
\newline
When no confusion is possible about the representation and the anti-representation, we shall write
$\mathcal{H}\otimes_{\psi}\mathcal{K}$ instead of
$\mathcal{H}\underset{\psi}{_\beta\otimes_\gamma}\mathcal{K}$, and $\xi\otimes_\psi\eta$ instead
of
$\xi\underset{\psi}{_\beta\otimes_\gamma}\eta$.
\newline
If $\theta\in Aut N$, then, using a remark made in \ref{spatial}, we get that the application which sends $\xi\underset{\psi}{_\beta\otimes_\gamma}\eta$ onto $\xi\underset{\psi\circ\theta}{_{\beta\circ\theta}\otimes_{\alpha\circ\theta}}\eta$ leads to a unitary from $\mathcal H\underset{\psi}{_\beta\otimes_\gamma}\mathcal K$ onto $\mathcal H\underset{\psi\circ\theta}{_{\beta\circ\theta}\otimes_{\alpha\circ\theta}}\mathcal K$. 
\newline
For any $\xi$ in $D(\mathcal{H}_\beta,
\psi^o)$, we define the bounded linear application $\lambda_\xi^{\beta, \gamma}$ from $\mathcal
K$ to
$\mathcal{H}\underset{\psi}{_\beta\otimes_\gamma}\mathcal{K}$ by, for all $\eta$ in $\mathcal K$,
$\lambda_\xi^{\beta, \gamma} (\eta)=\xi\underset{\psi}{_\beta\otimes_\gamma}\eta$. We shall write
$\lambda_\xi$ if no confusion is possible. We get ([EN], 3.10) :
\[\lambda_\xi^{\beta, \gamma}=R^{\beta, \psi^o}(\xi)\otimes_\psi 1_\mathcal K\]
where we recall the canonical identification (as left $N$-modules) of 
$L^2(N)\otimes_\psi\mathcal K$ with $\mathcal K$. We have :
\[(\lambda_\xi^{\beta, \gamma})^*\lambda_\xi^{\beta, \gamma}=\gamma(<\xi, \xi>_{\beta, \psi^o})\]
\newline
In ([S1] 2.1), the relative tensor product
$\mathcal{H}\underset{\psi}{_\beta\otimes_\gamma}\mathcal{K}$ is defined also, if
$\xi_1$, $\xi_2$ are in $\mathcal{H}$, $\eta_1$, $\eta_2$ are in $D(_\gamma\mathcal{K},\psi)$, by the
following formula :
\[(\xi_1\odot\eta_1 |\xi_2\odot\eta_2 )= (\beta(<\eta_1, \eta_2>_{\gamma,\psi})\xi_1 |\xi_2)\]
which leads to the the definition of a relative flip $\sigma_\psi$ which will be an isomorphism from
$\mathcal{H}\underset{\psi}{_\beta\otimes_\gamma}\mathcal{K}$ onto
$\mathcal{K}\underset{\psi^o}{_\gamma\otimes _\beta}\mathcal{H}$, defined, for any 
$\xi$ in $D(\mathcal{H}_\beta ,\psi^o )$, $\eta$ in $D(_\gamma \mathcal{K},\psi)$, by :
\[\sigma_\psi (\xi\otimes_\psi\eta)=\eta\otimes_{\psi^o}\xi\]
This allows us to define a relative flip $\varsigma_\psi$ from
$\mathcal{L}(\mathcal{H}\underset{\psi}{_\beta\otimes_\gamma}\mathcal{K})$ to $\mathcal{L}
(\mathcal{K}\underset{\psi^o}{_\gamma\otimes_\beta}\mathcal{H})$ which sends
$X$ in 
$\mathcal{L}(\mathcal{H}\underset{\psi}{_\beta\otimes_\gamma}\mathcal{K})$ onto
$\varsigma_\psi(X)=\sigma_\psi X\sigma_\psi^*$. Starting from another faithful semi-finite normal
weight $\psi'$, we get a von Neumann algebra
$\mathcal{L}(\mathcal{H}\underset{\psi'}{_\beta\otimes_\gamma}\mathcal{K})$ which is isomorphic to
$\mathcal{L}(\mathcal{H}\underset{\psi}{_\beta\otimes_\gamma}\mathcal{K})$, and a von Neumann
algebra $\mathcal{L} (\mathcal{K}\underset{\psi^{'o}}{_\gamma\otimes_\beta}\mathcal{H})$ which is
isomorphic to
$\mathcal{L} (\mathcal{K}\underset{\psi^o}{_\gamma\otimes_\beta}\mathcal{H})$; as we get that :
 \[\sigma_{\psi'}\circ U^{\psi, \psi'}_{\beta, \gamma}=U^{\psi^o, \psi'^o}_{\gamma, \beta}\]
 we see that these isomorphisms exchange $\varsigma_\psi$ and
$\varsigma_{\psi '}$. Therefore, the homomorphism $\varsigma_{\psi}$ can be denoted $\varsigma_N$
without any reference to a specific weight.
\newline
We may define, for any $\eta$ in $D(_\gamma\mathcal{K}, \psi)$, an application
$\rho_\eta^{\beta, \gamma}$ from $\mathcal H$ to
$\mathcal{H}\underset{\psi}{_\beta\otimes_\gamma}\mathcal{K}$ by
$\rho_\eta^{\beta, \gamma} (\xi)=\xi\underset{\psi}{_\beta\otimes_\gamma}\eta$. We shall write
$\rho_\eta$ if no confusion is possible. We get that :
\[(\rho_\eta^{\beta, \gamma})^*\rho_\eta^{\beta, \gamma}=\beta(<\eta, \eta>_{\gamma, \psi})\]
We recall, following
([S2], 2.2b) that, for all $\xi$ in $\mathcal{H}$, $\eta$ in $D(_\gamma\mathcal{K},\psi)$, $y$ in
$N$, analytic with respect to $\psi$, we have :
\[\beta (y)\xi \otimes_\psi\eta=\xi\otimes_\psi\gamma(\sigma^\psi_{-i/2}(y))\eta\]
Let $x$ be an element
of
$\mathcal{L}(\mathcal{H})$, commuting with the right action of $N$ on $\mathcal{H}_\beta$ (i.e. $x\in\beta(N)'$). It
is possible to define an operator $x\underset{\psi}{_\beta\otimes_\gamma} 1_{\mathcal{K}}$ on
$\mathcal{H}\underset{\psi}{_\beta\otimes_\gamma}
\mathcal{K}$. We can easily evaluate $\|x\underset{\psi}{_\beta\otimes_\gamma} 1_{\mathcal{K}}\|$, for any finite $J\subset I$, for any $\eta_i$ in $\mathcal K$, we have :
 \begin{multline*}
 ((x^*x\underset{\psi}{_\beta\otimes_\gamma} 1_{\mathcal{K}})(\Sigma_{i\in J}e_i\underset{\psi}{_\beta\otimes_\gamma}\eta_i)|(\Sigma_{i\in J}e_i\underset{\psi}{_\beta\otimes_\gamma}\eta_i))=\\
 =\Sigma_{i\in J}(\gamma(<xe_i, xe_i>_{\beta, \psi^o})\eta_i|\eta_i)\\
 \leq \|x\|^2\Sigma_{i\in J}(\gamma(<e_i, e_i>_{\beta, \psi^o})\eta_i|\eta_i)\
 =\|x\|^2\|\Sigma_{i\in J}e_i\underset{\psi}{_\beta\otimes_\gamma}\eta_i\|
 \end{multline*}
 from which we get $\|x\underset{\psi}{_\beta\otimes_\gamma} 1_{\mathcal{K}}\|\leq\|x\|$. 
 \newline
 By the same way, if $y$ commutes with the left action of $N$ on
$_\gamma\mathcal{K}$ (i.e. $y$ is in  $\gamma(N)'$), it is possible to define
$1_{\mathcal{H}}\underset{\psi}{_\beta\otimes_\gamma}y$ on
$\mathcal{H}\underset{\psi}{_\beta\otimes_\gamma} \mathcal{K}$, and by composition, it is possible
to define then
$x\underset{\psi}{_\beta\otimes_\gamma} y$. If we start from another faithful semi-finite normal
weight $\psi '$, the canonical isomorphism $U^{\psi, \psi'}_{\beta, \gamma}$ from $\mathcal{H}\underset{\psi}{_\beta\otimes_\gamma}
\mathcal{K}$ to $\mathcal{H}\underset{\psi'}{_\beta\otimes_\gamma} \mathcal{K}$ sends
$x\underset{\psi}{_\beta\otimes_\gamma} y$ on $x\underset{\psi'}{_\beta\otimes_\gamma} y$ ([S2],
2.3 and 2.6); therefore, this operator can be denoted $x\underset{N}{_\beta\otimes_\gamma} y$
without any reference to a specific weight, and we get $\|x\underset{N}{_\beta\otimes_\gamma} y\|\leq\|x\|\|y\|$. 
 \newline
 If $\theta\in Aut N$, the unitary from $\mathcal H\underset{\psi}{_\beta\otimes_\gamma}\mathcal K$ onto $\mathcal H\underset{\psi\circ\theta}{_{\beta\circ\theta}\otimes_{\alpha\circ\theta}}\mathcal K$ sends $x\underset{N}{_\beta\otimes_\gamma}y$ on $x\underset{N}{_{\beta\circ\theta}\otimes_{\gamma\circ\theta}}y$. 
 \newline
  With the notations of \ref{spatial}, let $(e_i)_{i\in I}$ a $(\beta, \psi^o)$-orthogonal basis of $\mathcal H$; 
 let us remark that, for all $\eta$ in $\mathcal K$, we have :
 \[e_i\underset{\psi}{_\beta\otimes_\gamma}\eta=e_i\underset{\psi}{_\beta\otimes_\gamma}\gamma(<e_i, e_i>_{\beta, \psi^o})\eta\]
On the other hand, $\theta^{\beta, \psi^o}(e_i, e_i)$ is an orthogonal projection, and so is $\theta^{\beta, \psi^o}(e_i, e_i)\underset{N}{_\beta\otimes_\gamma}1$; this last operator is the projection on the subspace $e_i\underset{\psi}{_\beta\otimes_\gamma}\gamma(<e_i, e_i>_{\beta, \psi^o})\mathcal K$ ([E2], 2.3) and, therefore, we get that $\mathcal H\underset{\psi}{_\beta\otimes_\gamma}\mathcal K$ is the orthogonal sum of the subspaces $e_i\underset{\psi}{_\beta\otimes_\gamma}\gamma(<e_i, e_i>_{\beta, \psi^o})\mathcal K$;  for any $\Xi$ in $\mathcal{H}\underset{\psi}{_\beta\otimes_\gamma}
\mathcal{K}$, there exist $\xi_i$ in $\mathcal K$, such that $\gamma(<e_i, e_i>_{\beta, \psi^o})\xi_i=\xi_i$ and $\Xi=\sum_i e_i\underset{\psi}{_\beta\otimes_\gamma}\xi_i$, from which we get that $\sum_i\|\xi_i\|^2=\|\Xi\|^2$. 
 \newline
Let us suppose now that $\mathcal{K}$ is a $N-P$ bimodule; that means that there exists a von
Neumann algebra $P$, and a non-degenerate normal anti-representation $\epsilon$ of $P$ on
$\mathcal{K}$, such that
$\epsilon (P)\subset\gamma (N)'$. We shall write then $_\gamma\mathcal{K}_\epsilon$. If $y$ is in $P$, we
have seen that it is possible to define then the operator
$1_{\mathcal{H}}\underset{\psi}{_\beta\otimes_\gamma}\epsilon (y)$ on
$\mathcal{H}\underset{\psi}{_\beta\otimes_\gamma}\mathcal{K}$, and we define this way a
non-degenerate normal antirepresentation of $P$ on
$\mathcal{H}\underset{\psi}{_\beta\otimes_\gamma}\mathcal{K}$, we shall call again $\epsilon$ for
simplification. If $\mathcal H$ is a $Q-N$ bimodule, then
$\mathcal{H}\underset{\psi}{_\beta\otimes_\gamma}\mathcal{K}$ becomes a $Q-P$ bimodule (Connes'
fusion of bimodules).
\newline
Taking a faithful semi-finite normal weight
$\nu$  on $P$, and a left $P$-module $_{\zeta}\mathcal{L}$ (i.e. a Hilbert space $\mathcal{L}$ and a normal
non-degenerate representation $\zeta$ of $P$ on $\mathcal{L}$), it is possible then to define
$(\mathcal{H}\underset{\psi}{_\beta\otimes_\gamma}\mathcal{K})\underset{\nu}{_\epsilon\otimes_\zeta}\mathcal{L}$.
Of course, it is possible also to consider the Hilbert space
$\mathcal{H}\underset{\psi}{_\beta\otimes_\gamma}(\mathcal{K}\underset{\nu}{_\epsilon\otimes_\zeta}\mathcal{L})$.
It can be shown that these two Hilbert spaces are isomorphics as $\beta (N)'-\zeta
(P)^{'o}$-bimodules. (In ([V1] 2.1.3), the proof, given for $N=P$ abelian can be used, without
modification, in that wider hypothesis). We shall write then
$\mathcal{H}\underset{\psi}{_\beta\otimes_\gamma}\mathcal{K}\underset{\nu}{_\epsilon\otimes_\zeta}\mathcal{L}$
without parenthesis, to emphazise this coassociativity property of the relative tensor
product.
\newline
Dealing now with that Hilbert space $\mathcal{H}\underset{\psi}{_\beta\otimes_\gamma}\mathcal{K}\underset{\nu}{_\epsilon\otimes_\zeta}\mathcal{L}$, there exist different flips, and it is necessary to be careful with notations. For instance, $1\underset{\psi}{_\beta\otimes\sigma_\nu}$ is the flip from this Hilbert space onto $\mathcal{H}\underset{\psi}{_\beta\otimes_\gamma}(\mathcal L\underset{\nu^o}{_\zeta\otimes_\epsilon}\mathcal K)$, where $\gamma$ is here acting on the second leg of $\mathcal L\underset{\nu^o}{_\zeta\otimes_\epsilon}\mathcal K$ (and should therefore be written $1\underset{\nu^o}{_\zeta\otimes_\epsilon}\gamma$, but this will not be done for obvious reasons). Here, the parenthesis remains, because there is no associativity rule, and to remind that $\gamma$ is not acting on $\mathcal L$. The adjoint of $1\underset{\psi}{_\beta\otimes\sigma_\nu}$ is $1\underset{\psi}{_\beta\otimes\sigma_{\nu^o}}$. 
\newline
The same way, we can consider $\sigma_\psi\underset{\nu}{_\epsilon\otimes_\zeta}1$ from $\mathcal{H}\underset{\psi}{_\beta\otimes_\gamma}\mathcal{K}\underset{\nu}{_\epsilon\otimes_\zeta}\mathcal{L}$ onto $(\mathcal K\underset{\psi^o}{_\gamma\otimes_\beta}\mathcal H)\underset{\nu}{_\epsilon\otimes_\zeta}\mathcal L$. 
 \newline
 Another kind of flip sends $\mathcal{H}\underset{\psi}{_\beta\otimes_\gamma}(\mathcal L\underset{\nu^o}{_\zeta\otimes_\epsilon}\mathcal K)$ onto $\mathcal L\underset{\nu^o}{_\zeta\otimes_\epsilon}(\mathcal{H}\underset{\psi}{_\beta\otimes_\gamma}\mathcal{K})$. We shall denote this application $\sigma^{1,2}_{\gamma, \epsilon}$ (and its adjoint $\sigma^{1,2}_{\epsilon, \gamma}$), in order to emphasize that we are exchanging the first and the second leg, and the representation $\gamma$ and $\epsilon$ on the third leg. 
 \newline
 If $\pi$ denotes the canonical left representation of $N$ on the Hilbert space $L^2(N)$, then it is straightforward to verify that the application which sends, for all $\xi$ in $\mathcal {H}$, $\chi$ normal faithful semi-finite weight on $N$, and $x$ in $\gN_\chi$, the vector $\xi{}_\beta\underset{\chi}{\otimes}{}_\pi J_\chi\Lambda_\chi (x)$ on $\beta(x^*)\xi$, gives an isomorphism of $\mathcal{H}{}_\beta\underset{\chi}{\otimes}{}_\pi L^2(N)$ on $\mathcal{H}$, which will send the antirepresentation of $N$ given by $n\mapsto 1_{\mathcal{H}}{}_\beta\underset{\chi}{\otimes}{}_\pi J_\chi n^*J_\chi$ on $\beta$
 \newline
 If $\mathcal K$ is a Hilbert space on which there exists a non-degenerate representation $\gamma$ of $N$, then $\mathcal K$ is a $N-\gamma(N)'^o$ bimodule, and the conjugate Hilbert space $\overline{\mathcal K}$ is a $\gamma(N)'-N^o$ bimodule, and, ([S2]), for any normal faithful semi-finite weight $\phi$ on $\gamma(N)'$, the fusion $_\gamma\mathcal K\underset{\phi^o}{\otimes}\overline{\mathcal K}_\gamma$ is isomorphic to the standard space $L^2(N)$, equipped with its standard left and right representation. 
 \newline
 Using that remark, one gets for any $x\in \beta(N)'$ :
 \[\|x\underset{N}{_\beta\otimes_\gamma}1_{\mathcal K}\|\leq\|x\underset{N}{_\beta\otimes_\gamma}1_{\mathcal K}\underset{\gamma(N)'^o}{\otimes}1_{\overline{\mathcal K}}\|=\|x\underset{N}{_\beta\otimes}1_{L^2(N)}\|=\|x\|\]
from which we have $\|x\underset{N}{_\beta\otimes_\gamma}1_{\mathcal K}\|=\|x\|$.

\subsection{Fiber product [V1], [EV]} 
\label{fiber}
Let us follow the notations of \ref{rel}; let now
$M_1$ be a von Neumann algebra on $\mathcal{H}$, such that $\beta (N)\subset
M_1$, and $M_2$ be a von Neumann algebra on $\mathcal{K}$, such that $\gamma (N)\subset
M_2$. The von Neumann algebra generated by all elements $x\underset{N}{_\beta\otimes_\gamma} y$,
where
$x$ belongs to $M'_1$, and $y$ belongs $M'_2$ will be denoted
$M'_1\underset{N}{_\beta\otimes_\gamma} M'_2$ (or $M'_1\otimes_N M'_2$ if no confusion if
possible), and will be called the relative tensor product of
$M'_1$ and $M'_2$ over $N$. The commutant of this algebra will be denoted 
$M_1\underset{N}{_\beta *_\gamma} M_2$ (or $M_1*_N M_2$ if no confusion is possible) and called the
fiber product of $M_1$ and
$M_2$, over
$N$. If $\theta\in Aut N$, using a remark made in \ref{rel}, we get that the von Neumann algebras $M_1\underset{N}{_\beta*_\gamma}M_2$ and $M_1\underset{N}{_{\beta\circ\theta}*_{\gamma\circ\theta}}M_2$ are spatially isomorphic, and we shall identify them. 
\newline
It is straightforward to verify that, if $P_1$ and $P_2$ are two other von Neumann
algebras satisfying the same relations with $N$, we have 
\[M_1*_N M_2\cap P_1*_N P_2=(M_1\cap P_1)*_N (M_2\cap P_2)\]
Moreover, we get that $\varsigma_N (M_1\underset{N}{_\beta *_\gamma}
M_2)=M_2\underset{N^o}{_\gamma *_\beta}M_1$.
\newline
In particular, we have :
\[(M_1\cap \beta (N)')\underset{N}{_\beta\otimes_\gamma} (M_2\cap \gamma (N)')\subset
M_1\underset{N}{_\beta *_\gamma} M_2\] and :
\[M_1\underset{N}{_\beta *_\gamma} \gamma(N)=(M_1\cap\beta (N)')\underset{N}{_\beta\otimes_\gamma} 1\]
More generally, if
$\beta$ is a non-degenerate normal involutive antihomomorphism from
$N$ into a von Neumann algebra
$M_1$, and
$\gamma$ a non-degenerate normal involutive homomorphism from $N$ into a von Neumann
algebra
$M_2$, it is possible
to define, without any reference to a specific Hilbert space, a von Neumann algebra
$M_1\underset{N}{_\beta *_ \gamma}M_2$. 
\newline
Moreover, if now $\beta '$ is a non-degenerate normal involutive antihomomorphism from $N$ into
another von Neumann algebra
$P_1$,
$\gamma '$ a non-degenerate normal involutive homomorphism from $N$ into another
von Neumann algebra $P_2$, $\Phi$ a normal involutive homomorphism from $M_1$ into $P_1$ such that
$\Phi\circ\beta =\beta '$, and $\Psi$ a normal involutive homomorphism from $M_2$ into $P_2$ such that
$\Psi\circ\gamma=\gamma'$, it is possible then to define a normal involutive homomorphism (the proof
given in ([S1] 1.2.4) in the case when $N$ is abelian can be extended without modification in the
general case) :
\[\Phi\underset{N}{_\beta *_\gamma}\Psi 
: M_1\underset{N}{_\beta
*_\gamma}M_2\mapsto P_1\underset{N}{_{\beta '}*_{\gamma '}}P_2\]
Let $\Phi$ be in $Aut M_1$, $\Psi$ in $Aut M_2$, and let $\theta\in Aut N$ be such that $\Phi\circ\beta=\beta\circ\theta$ and $\Psi\circ\gamma=\gamma\circ\theta$, then, using the identification between $M_1\underset{N}{_{\beta }*_{\gamma }}M_2$ and $M_1\underset{N}{_{\beta \circ\theta}*_{\gamma \circ\theta}}M_2$, we get the existence of an automorphism $\Phi\underset{N}{_\beta *_\gamma}\Psi $ of $M_1\underset{N}{_\beta
*_\gamma}M_2$. 
\newline
In the case when $_\gamma\mathcal{K}_\epsilon$ is a $N-P^o$ bimodule as explained in \ref{rel} and
$_\zeta\mathcal{L}$ a $P$-module, if
$\gamma (N)\subset M_2$ and $\epsilon (P)\subset M_2$, and if $\zeta (P)\subset M_3$, where $M_3$ is
a von Neumann algebra on $\mathcal{L}$, it is possible to consider then $(M_1\underset{N}{_\beta
*_\gamma}M_2)\underset{P}{_\epsilon *_\zeta}M_3$ and $M_1\underset{N}{_\beta
*_\gamma}(M_2\underset{P}{_\epsilon *_\zeta}M_3)$. The coassociativity property for relative tensor
products leads then to the isomorphism of these von Neumann algebra we shall write now 
$M_1\underset{N}{_\beta
*_\gamma}M_2\underset{P}{_\epsilon *_\zeta}M_3$ without parenthesis.

\subsection{Slice maps [E3]}
\label{slice}
Let $A$ be in $M_1\underset{N}{_\beta *_\gamma}M_2$, $\psi$ a normal faithful semi-finite weight on $N$, $\mathcal{H}$ an Hilbert space on which $M_1$ is acting, $\mathcal{K}$ an Hilbert space on which $M_2$ is acting, and let $\xi_1$, $\xi_2$ be in
$D(\mathcal{H}_\beta,\psi^o)$; let us define :
\[(\omega_{\xi_1, \xi_2}\underset{\psi}{_\beta*_\gamma}id)(A)=(\lambda^{\beta, \gamma}_{\xi_2})^*A\lambda^{\beta, \gamma}_{\xi_1}\]
We define this way $(\omega_{\xi_1, \xi_2}\underset{\psi}{_\beta*_\gamma}id)(A)$ as a bounded operator on $\mathcal{K}$,
which belongs to $M_2$, such that :
\[((\omega_{\xi_1, \xi_2}\underset{\psi}{_\beta*_\gamma}id)(A)\eta_1|\eta_2)=
(A(\xi_1\underset{\psi}{_\beta\otimes_\gamma}\eta_1)|
\xi_2\underset{\psi}{_\beta\otimes_\gamma}\eta_2)\]
One should note that $(\omega_{\xi_1, \xi_2}\underset{\psi}{_\beta*_\gamma}id)(1)=\gamma (<\xi_1, \xi_2 >_{\beta, \psi^o})$. 
\newline
Let us define the same way, for any $\eta_1$, $\eta_2$ in
$D(_\gamma\mathcal{K}, \psi)$:
\[(id\underset{\psi}{_\beta*_\gamma}\omega_{\eta_1, \eta_2})(A)=(\rho^{\beta, \gamma}_{\eta_2})^*A\rho^{\beta, \gamma}_{\eta_1}\]
which belongs to $M_1$. 
\newline
We therefore have a Fubini formula for these slice maps : for any $\xi_1$, $\xi_2$ in
$D(\mathcal{H}_\beta,\psi^o)$, $\eta_1$, $\eta_2$ in $D(_\gamma\mathcal{K}, \psi)$, we have :
\[<(\omega_{\xi_1, \xi_2}\underset{\psi}{_\beta*_\gamma}id)(A), \omega_{\eta_1, \eta_2}>=<(id\underset{\psi}{_\beta*_\gamma}\omega_{\eta_1,
\eta_2})(A),\omega_{\xi_1, \xi_2}>\]
Let $\phi_1$ be a normal semi-finite weight on
$M_1^+$, and $A$ be a positive element of the fiber product
$M_1\underset{N}{_\beta*_\gamma}M_2$, then we may define an element of the extended positive part
of $M_2$, denoted
$(\phi_1\underset{\psi}{_\beta*_\gamma}id)(A)$, such that, for all $\eta$ in $D(_\gamma L^2(M_2), \psi)$, we have :
\[\|(\phi_1\underset{\psi}{_\beta*_\gamma}id)(A)^{1/2}\eta\|^2=\phi_1(id\underset{\psi}{_\beta*_\gamma}\omega_\eta)(A)\]
Moreover, then, if $\phi_2$ is a normal semi-finite weight on $M_2^+$, we have :
\[\phi_2(\phi_1\underset{\psi}{_\beta*_\gamma}id)(A)=\phi_1(id\underset{\psi}{_\beta*_\gamma}\phi_2)(A)\]
and if $\omega_i$ are in $M_{1*}$ such that $\phi_1=sup_i\omega_i$, we
have $(\phi_1\underset{\psi}{_\beta*_\gamma}id)(A)=sup_i(\omega_i\underset{\psi}{_\beta*_\gamma}id)(A)$.
\newline
Let now $P_1$ be a von Neuman algebra such that :
\[\beta(N)\subset P_1\subset M_1\]
and let $\Phi_i$ ($i=1,2$)
be a normal faithful semi-finite operator valued weight from $M_i$ to $P_i$; for any positive
operator $A$ in the fiber product
$M_1\underset{N}{_\beta*_\gamma}M_2$, there exists an element $(\Phi_1\underset{\psi}{_\beta*_\gamma}id)(A)$
of the extended positive part
of $P_1\underset{N}{_\beta*_\gamma}M_2$, such that ([E3], 3.5), for all
$\eta$ in $D(_\gamma L^2(M_2), \psi)$, and $\xi$ in $D(L^2(P_1)_\beta, \psi^o)$, we have :
\[\|(\Phi_1\underset{\psi}{_\beta*_\gamma}id)(A)^{1/2}(\xi\underset{\psi}{_\beta\otimes_\gamma}\eta)\|^2=
\|\Phi_1(id\underset{\psi}{_\beta*_\gamma}\omega_\eta)(A)^{1/2}\xi\|^2\]
If $\phi$ is a normal semi-finite weight on $P$, we have :
\[(\phi\circ\Phi_1\underset{\psi}{_\beta*_\gamma}id)(A)=(\phi\underset{\psi}{_\beta*_\gamma}id)
(\Phi_1\underset{\psi}{_\beta*_\gamma}id)(A)\]
We define the same way an element $(id\underset{\psi}{_\beta*_\gamma}\Phi_2)(A)$ of the extended positive part
of
$M_1\underset{N}{_\gamma*_\beta}P_2$, and we have :
\[(id\underset{\psi}{_\beta*_\gamma}\Phi_2)((\Phi_1\underset{\psi}{_\beta*_\gamma}id)(A))=
(\Phi_1\underset{\psi}{_\beta*_\gamma}id)((id\underset{\psi}{_\beta*_\gamma}\Phi_2)(A))\]
 Considering now an element $x$ of $M_1{}_\beta\underset{\psi}{*}{}_\pi \pi(N)$, which can be identified 
(\ref{fiber}) to $M_1\cap\beta(N)'$, we get that, for $e$ in $\gN_\psi$, we have  \[(id_\beta\underset{\psi}{*}{}_\pi\omega_{J_\psi \Lambda_{\psi}(e)})(x)=\beta(ee^*)x\]
 Therefore, by increasing limits, we get that $(id_\beta\underset{\psi}{*}{}_\pi\psi)$ is the injection of $M_1\cap\beta(N)'$ into $M_1$.  More precisely, if $x$ belongs to $M_1\cap\beta(N)'$, we have :
 \[(id_\beta\underset{\psi}{*}{}_\pi\psi)(x{}_\beta\underset{\psi}{\otimes}{}_\pi 1)=x\]
 \newline
 Therefore, if $\Phi_2$ is a normal faithful semi-finite operator-valued weight from $M_2$ onto $\gamma(N)$, we get that, for all $A$ positive in $M_1\underset{N}{_\beta*_\gamma}M_2$, we have :
 \[(id_\beta\underset{\psi}{*}{}_\gamma\psi\circ\Phi_2)(A){}_\beta\underset{\psi}{\otimes}{}_\gamma 1=
 (id_\beta\underset{\psi}{*}{}_\gamma\Phi_2)(A)\]

With the notations of \ref{spatial}, let $(e_i)_{i\in I}$ be a $(\beta, \psi^o)$-orthogonal basis of $\mathcal H$; using the fact (\ref{rel}) that, for all $\eta$ in $\mathcal K$, we have :\[e_i\underset{\psi}{_\beta\otimes_\gamma}\eta=e_i\underset{\psi}{_\beta\otimes_\gamma}\gamma(<e_i, e_i>_{\beta, \psi^o})\eta\]
we get that, for all $X$ in $M_1\underset{N}{_\beta*_\gamma}M_2$, $\xi$ in $D(\mathcal H_\beta, \psi^o)$, we have 
\[(\omega_{\xi, e_i}\underset{\psi}{_\beta *_\gamma}id)(X)=\gamma(<e_i, e_i>_{\beta, \psi^o})(\omega_{\xi, e_i}\underset{\psi}{_\beta *_\gamma}id)(X)\]

\subsection{Vaes' Radon-Nikodym theorem}
\label{Vaes}
In [V] is proved a very nice Radon-Nikodym theorem for two normal faithful semi-finite weights on a von Neumann algebra $M$. If $\Phi$ and $\Psi$ are such weights, then are equivalent :
\newline
- the two modular automorphism groups $\sigma^\Phi$ and $\sigma^\Psi$ commute;
\newline
- the Connes' derivative $[D\Psi : D\Phi]_t$ is of the form :
\[[D\Psi : D\Phi]_t=\lambda^{it^2/2}\delta^{it}\]
where $\lambda$ is a non-singular positive operator affiliated to $Z(M)$, and $\delta$ is a non-singular positive operator affiliated to $M$. 
\newline
It is then easy to verify that $\sigma^\Phi_t(\delta^{is})=\lambda^{ist}\delta^{is}$, and that 
\[[D\Phi\circ\sigma^\Psi_t:D\Phi]_s=\lambda^{ist}\]
\[[D\Psi\circ\sigma^\Phi_t:D\Psi]_s=\lambda^{-ist}\]
Moreover, we have also, for any $x\in M^+$ :
\[\Psi(x)=lim_n\Phi((\delta^{1/2}e_n)x(\delta^{1/2}e_n))\]
where the $e_n$ are self-adjoint elements of $M$ given by the formula :
\[e_n=a_n\int_{\mathbb{R}^2}e^{-n^2x^2-n^4y^4}\lambda^{ix}\delta^{iy}dxdy\]
where we put $a_n=2n^2\Gamma(1/2)^{-1}\Gamma(1/4)^{-1}$. 
The operators $e_n$ are analytic with respect to $\sigma^\Phi$ and such that, for any $z\in \mathbb{C}$,  the sequence $\sigma_z^\Phi(e_n)$ is bounded and strongly converges to $1$. 
\newline
In that situation, we shall write $\Psi=\Phi_\delta$ and call $\delta$ the modulus of $\Psi$ with respect to $\Phi$; $\lambda$ will be called the scaling operator of $\Psi$ with respect to $\Phi$. 
\newline
Moreover, if $a\in M$ is such that $a\delta^{1/2}$ is bounded and its closure $\overline{a\delta^{1/2}}$ belongs to $\gN_\Phi$, then $a$ belongs to $\gN_\Psi$. We may then identify $\Lambda_\Psi(a)$ with $\Lambda_\Phi(\overline{a\delta^{1/2}})$, $J_\Psi$ with $\lambda^{i/4}J_\Phi$, $\Delta_\Psi$ with $\overline{J_\Phi\delta^{-1}J_\Phi\delta\Delta_\Phi}$.


\section{Hopf-bimodules and Pseudo-multiplicative unitary}
\label{pmu}
In this chapter, we recall the definition of Hopf-bimodules
(\ref{Hbimod}), the definition of a pseudo-multiplicative unitary (\ref{defmult}), give the fundamental example given by groupo\"{\i}ds (\ref{gd}), and construct the algebras and the Hopf-bimodules "generated by the left (resp. right) leg" of a pseudo-multiplicative unitary (\ref{AW}). We recall the definition of left- (resp. right-) invariant operator-valued weights on a Hopf-bimodule; if we have both operator-valued weights, we then recall Lesieur's construction of a pseudo-multiplicative unitary. 

\subsection{Definition}
\label{Hbimod}
A quintuplet $(N, M, \alpha, \beta, \Gamma)$ will be called a Hopf-bimodule, following ([Val1], [EV] 6.5), if
$N$,
$M$ are von Neumann algebras, $\alpha$ a faithful non-degenerate representation of $N$ into $M$, $\beta$ a
faithful non-degenerate anti-representation of
$N$ into $M$, with commuting ranges, and $\Gamma$ an injective involutive homomorphism from $M$
into
$M\underset{N}{_\beta *_\alpha}M$ such that, for all $X$ in $N$ :
\newline
(i) $\Gamma (\beta(X))=1\underset{N}{_\beta\otimes_\alpha}\beta(X)$
\newline
(ii) $\Gamma (\alpha(X))=\alpha(X)\underset{N}{_\beta\otimes_\alpha}1$ 
\newline
(iii) $\Gamma$ satisfies the co-associativity relation :
\[(\Gamma \underset{N}{_\beta *_\alpha}id)\Gamma =(id \underset{N}{_\beta *_\alpha}\Gamma)\Gamma\]
This last formula makes sense, thanks to the two preceeding ones and
\ref{fiber}\vspace{5mm}.\newline
If $(N, M, \alpha, \beta, \Gamma)$ is a Hopf-bimodule, it is clear that
$(N^o, M, \beta, \alpha,
\varsigma_N\circ\Gamma)$ is another Hopf-bimodule, we shall call the symmetrized of the first
one. (Recall that $\varsigma_N\circ\Gamma$ is a homomorphism from $M$ to
$M\underset{N^o}{_r*_s}M$).
\newline
If $N$ is abelian, $\alpha=\beta$, $\Gamma=\varsigma_N\circ\Gamma$, then the quadruplet $(N, M, \alpha, \alpha,
\Gamma)$ is equal to its symmetrized Hopf-bimodule, and we shall say that it is a symmetric
Hopf-bimodule\vspace{5mm}.\newline
Let $\mathcal G$ be a groupo\"{\i}d, with $\mathcal G^{(0)}$ as its set of units, and let us denote
by $r$ and $s$ the range and source applications from $\mathcal G$ to $\mathcal G^{(0)}$, given by
$xx^{-1}=r(x)$ and $x^{-1}x=s(x)$. As usual, we shall denote by $\mathcal G^{(2)}$ (or $\mathcal
G^{(2)}_{s,r}$) the set of composable elements, i.e. 
\[\mathcal G^{(2)}=\{(x,y)\in \mathcal G^2; s(x)=r(y)\}\]
In [Y] and [Val1] were associated to a measured groupo\"{\i}d $\mathcal G$, equipped with a Haar system $(\lambda^u)_{u\in \mathcal G ^{(0)}}$ and a quasi-invariant measure $\mu$ on $\mathcal G ^{(0)}$ (see [R1],
[R2], [C2] II.5 and [AR] for more details, precise definitions and examples of groupo\"{\i}ds) two
Hopf-bimodules : 
\newline
The first one is $(L^\infty (\mathcal G^{(0)}, \mu), L^\infty (\mathcal G, \nu), r_{\mathcal G}, s_{\mathcal G}, \Gamma_{\mathcal
G})$, where $\nu$ is the measure constructed on $\mathcal G $ using $\mu$ and the Haar system $(\lambda^u)_{u\in \mathcal G ^{(0)}}$, where we define $r_{\mathcal G}$ and $s_{\mathcal G}$ by writing , for $g$ in $L^\infty (\mathcal G^{(0)})$ :
\[r_{\mathcal G}(g)=g\circ r\]
\[s_{\mathcal G}(g)=g\circ s\]
 and where
$\Gamma_{\mathcal G}(f)$, for $f$ in $L^\infty (\mathcal G)$, is the function defined on $\mathcal G^{(2)}$ by $(s,t)\mapsto f(st)$;
$\Gamma_{\mathcal G}$ is then an involutive homomorphism from $L^\infty (\mathcal G)$ into $L^\infty
(\mathcal G^2_{s,r})$ (which can be identified to
$L^\infty (\mathcal G){_s*_r}L^\infty (\mathcal G)$).
\newline
The second one is symmetric; it is $(L^\infty (\mathcal G^{(0)}), \mathcal L(\mathcal G), r_{\mathcal G}, r_{\mathcal G},
\widehat{\Gamma_{\mathcal G}})$, where
$\mathcal L(\mathcal G)$ is the von Neumann algebra generated by the convolution algebra associated to the
groupo\"{\i}d
$\mathcal G$, and $\widehat{\Gamma_{\mathcal G}}$ has been defined in [Y] and
[Val1]. 

\subsection{Definition}
\label{defmult}
Let $N$ be a von Neumann algebra; let
$\gH$ be a Hilbert space on which $N$ has a non-degenerate normal representation $\alpha$ and two
non-degenerate normal anti-representations $\hat{\beta}$ and $\beta$. These 3 applications
are supposed to be injective, and to commute two by two.  Let $\nu$ be a normal semi-finite faithful weight on
$N$; we can therefore construct the Hilbert spaces
$\gH\underset{\nu}{_\beta\otimes_\alpha}\gH$ and
$\gH\underset{\nu^o}{_\alpha\otimes_{\hat{\beta}}}\gH$. A unitary $W$ from
$\gH\underset{\nu}{_\beta\otimes_\alpha}\gH$ onto
$\gH\underset{\nu^o}{_\alpha\otimes_{\hat{\beta}}}\gH$.
will be called a pseudo-multiplicative unitary over the basis $N$, with respect to the
representation $\alpha$, and the anti-representations $\hat{\beta}$ and $\beta$ (we shall say it is an $(\alpha, \hat{\beta}, \beta)$-pseudo-multiplicative unitary), if :
\newline
(i) $W$ intertwines $\alpha$, $\hat{\beta}$, $\beta$  in the following way :
\[W(\alpha
(X)\underset{N}{_\beta\otimes_\alpha}1)=
(1\underset{N^o}{_\alpha\otimes_{\hat{\beta}}}\alpha(X))W\]
\[W(1\underset{N}{_\beta\otimes_\alpha}\beta
(X))=(1\underset{N^o}{_\alpha\otimes_{\hat{\beta}}}\beta (X))W\]
\[W(\hat{\beta}(X) \underset{N}{_\beta\otimes_\alpha}1)=
(\hat{\beta}(X)\underset{N^o}{_\alpha\otimes_{\hat{\beta}}}1)W\]
\[W(1\underset{N}{_\beta\otimes_\alpha}\hat{\beta}(X))=
(\beta(X)\underset{N^o}{_\alpha\otimes_{\hat{\beta}}}1)W\]
(ii) The operator $W$ satisfies :
\begin{multline*}
(1_\gH\underset{N^o}{_\alpha\otimes_{\hat{\beta}}}W)
(W\underset{N}{_\beta\otimes_\alpha}1_{\gH})=\\
=(W\underset{N^o}{_\alpha\otimes_{\hat{\beta}}}1_{\gH})
\sigma^{2,3}_{\alpha, \beta}(W\underset{N}{_{\hat{\beta}}\otimes_\alpha}1)
(1_{\gH}\underset{N}{_\beta\otimes_\alpha}\sigma_{\nu^o})
(1_{\gH}\underset{N}{_\beta\otimes_\alpha}W)
\end{multline*}
Here, $\sigma^{2,3}_{\alpha, \beta}$
goes from $(H\underset{\nu^o}{_\alpha\otimes_{\hat{\beta}}}H)\underset{\nu}{_\beta\otimes_\alpha}H$ to $(H\underset{\nu}{_\beta\otimes_\alpha}H)\underset{\nu^o}{_\alpha\otimes_{\hat{\beta}}}H$, 
and $1_{\gH}\underset{N}{_\beta\otimes_\alpha}\sigma_{\nu^o}$ goes from $H\underset{\nu}{_\beta\otimes_\alpha}(H\underset{\nu^o}{_\alpha\otimes_{\hat{\beta}}}H)$ to $H\underset{\nu}{_\beta\otimes_\alpha}H\underset{\nu}{_{\hat{\beta}}\otimes_\alpha}H$. 
\newline
All the properties supposed in (i) allow us to write such a formula, which will be called the
"pentagonal relation". 
\newline
One should note that this definition is different from the definition introduced in [EV] (and repeated afterwards). It is in fact the same formula, the new writing 
\[\sigma^{2,3}_{\alpha, \beta}
(W\underset{N}{_{\hat{\beta}}\otimes_\alpha}1)
(1_{\gH}\underset{N}{_\beta\otimes_\alpha}\sigma_{\nu^o})\]
is here replacing the rather akward writing 
\[(\sigma_{\nu^o}\underset{N^o}{_\alpha\otimes_{\hat{\beta}}}1_{\gH})
(1_{\gH}\underset{N^o}{_\alpha\otimes_{\hat{\beta}}}W)
\sigma_{2\nu}
(1_{\gH}\underset{N}{_\beta\otimes_\alpha}\sigma_{\nu^o})\]
but denotes the same operator, and we suggest the reader to convince himself of this easy fact. 
\newline
All the properties supposed in (i) allow us to write such a formula, which will be called the
"pentagonal relation". 
\newline
If we start from another normal semi-finite faithful weight $\nu'$ on $N$, we may define, using \ref{rel}, another unitary $W^{\nu'}=U^{\nu^o,
\nu^{'o}}_{\alpha, {\hat{\beta}}}WU^{\nu', \nu}_{\beta, \alpha}$ from $\gH\underset{\nu'}{_\beta\otimes_\alpha}\gH$ onto
$\gH\underset{\nu^{'o}}{_\alpha\otimes_{\hat{\beta}}}\gH$. The formulae which link these isomorphims between relative product Hilbert spaces and the
relative flips allow us to check that this operator $W^{\nu'}$ is also pseudo-multiplicative; which can be resumed in saying that a
pseudo-multiplicative unitary does not depend on the choice of the weight on $N$. 
\newline
If $W$ is an $(\alpha, \hat{\beta}, \beta)$-pseudo-multiplicative unitary, then the unitary $\sigma_\nu W^*\sigma_\nu$ from $\gH\underset{\nu}{_{\hat{\beta}}\otimes_\alpha}\gH$ to $\gH\underset{\nu^o}{_\alpha\otimes_\beta}\gH$ is an $(\alpha, \beta, \hat{\beta})$-pseudo-multiplicative unitary, called the dual of $W$.

\subsection{Algebras and Hopf-bimodules associated to a pseudo-multiplicative unitary}
\label{AW}
For $\xi_2$ in $D(_\alpha\gH, \nu)$, $\eta_2$ in $D(\gH_{\hat{\beta}}, \nu^o)$, the operator $(\rho_{\eta_2}^{\alpha,
\hat{\beta}})^*W\rho_{\xi_2}^{\beta, \alpha}$ will be written $(id*\omega_{\xi_2, \eta_2})(W)$; we have, therefore, for all
$\xi_1$, $\eta_1$ in $\gH$ :
\[((id*\omega_{\xi_2, \eta_2})(W)\xi_1|\eta_1)=(W(\xi_1\underset{\nu}{_\beta\otimes_\alpha}\xi_2)|
\eta_1\underset{\nu^o}{_\alpha\otimes_{\hat{\beta}}}\eta_2)\]
and, using the intertwining property of $W$ with $\hat{\beta}$, we easily get that $(id*\omega_{\xi_2, \eta_2})(W)$ belongs
to $\hat{\beta} (N)'$. 
\newline
If $x$ belongs to $N$, we have $(id*\omega_{\xi_2, \eta_2})(W)\alpha (x)=(id*\omega_{\xi_2, \alpha(x^*)\eta_2})(W)$, and $\beta(x)(id*\omega_{\xi_2, \eta_2})(W)=(id*\omega_{\hat{\beta}(x)\xi_2, \eta_2})(W)$. 
\newline
If $\xi$ belongs to $D(_\alpha\gH, \nu)\cap D(\gH_{\hat{\beta}}, \nu^o)$, we shall write $(id*\omega_\xi)(W)$ instead of
$(id*\omega_{\xi, \xi})(W)$.
\newline
We shall write $A_w(W)$ the weak closure of the linear span of these operators, which are right $\alpha(N)$-modules and left $\beta(N)$-modules. Applying ([E2] 3.6), we get that $A_w(W)^*$ and $A_w(W)$ are non-degenerate algebras (one should note that the notations of ([E2]) had been changed in order to fit with Lesieur's notations). 
We shall write $\mathcal A(W)$ the von Neumann algebra generated by $A_w(W)$ .
We then have $\mathcal A(W)\subset\hat{\beta}(N)'$.
\newline
For $\xi_1$ in $D(\gH_\beta,\nu^o)$, $\eta_1$ in $D(_\alpha\gH, \nu)$, the operator $(\lambda_{\eta_1}^{\alpha,
\hat{\beta}})^*W\lambda_{\xi_1}^{\beta, \alpha}$ will be written $(\omega_{\xi_1,
\eta_1}*id)(W)$; we have,
therefore, for all
$\xi_2$,
$\eta_2$ in
$\gH$ :
\[((\omega_{\xi_1,\eta_1}*id)(W)\xi_2|\eta_2)=(W(\xi_1\underset{\nu}{_\beta\otimes_\alpha}\xi_2)|
\eta_1\underset{\nu^o}{_\alpha\otimes_{\hat{\beta}}}\eta_2)\]
and, using the intertwining property of $W$ with $\beta$, we easily get that $(\omega_{\xi_1,
\eta_1}*id)(W)$ belongs to $\beta(N)'$. If $\xi$ belongs to $D(\gH_\beta, \nu^o)\cap D(_\alpha\gH, \nu)$, we shall write
$(\omega_\xi*id)(W)$ instead of $(\omega_{\xi, \xi}*id)(W)$. 
\newline
We shall write $\widehat{A_w(W)}$ the weak closure of the linear span of these operators. It is clear that this weakly closed subspace is a non degenarate algebra; following ([EV] 6.1 and 6.5), we shall write $\widehat{\mathcal A(W})$ the von Neumann algebra generated by  $\widehat{A_w(W)}$. We then have $\widehat{\mathcal A(W)}\subset\beta(N)'$. 
\newline
In ([EV] 6.3 and 6.5), using the pentagonal equation, we got
that
$(N,\mathcal A(W),\alpha,\beta,\Gamma)$, and
$(N^o,\widehat{\mathcal A(W)}, \hat{\beta}, \alpha, \widehat{\Gamma})$ are Hopf-bimodules, where $\Gamma$ and
$\widehat{\Gamma}$ are defined, for any $x$ in $\mathcal A(W)$ and $y$ in $\widehat{\mathcal
A(W)}$, by :
\[\Gamma(x)=W^*(1\underset{N^o}{_\alpha\otimes_{\hat{\beta}}}x)W\]
\[\widehat{\Gamma}(y)=W(y\underset{N}{_\beta\otimes_\alpha}1)W^*\]
In ([EV] 6.1(iv)), we had obtained that $x$ in $\mathcal L(\gH)$ belongs to $\mathcal A(W)'$ if and only if $x$ belongs to $\alpha(N)'\cap
\beta(N)'$ and verifies \[(x\underset{N^o}{_\alpha\otimes_{\hat{\beta}}}1)W=W(x\underset{N}{_\beta\otimes_\alpha}1)\]
We obtain the same way
that $y$ in $\mathcal L(\gH)$ belongs to $\widehat{\mathcal A(W)}'$ if and only if $y$ belongs to $\alpha(N)'\cap
\hat{\beta}(N)'$ and verify $(1\underset{N^o}{_\alpha\otimes_{\hat{\beta}}}y)W=W(1\underset{N}{_\beta\otimes_\alpha}y)$.  \newline
Moreover, we get that $\alpha(N)\subset\mathcal A\cap\widehat{\mathcal A}$, $\beta(N)\subset\mathcal A$,
$\hat{\beta}(N)\subset\widehat{\mathcal A}$, and, for all $x$ in $N$ :
\[\Gamma (\alpha (x))=\alpha (x)\underset{N}{_\beta\otimes_\alpha}1\]
\[\Gamma (\beta (x))=1\underset{N}{_\beta\otimes_\alpha}\beta (x)\]
\[\widehat{\Gamma}(\alpha(x))=1\underset{N^o}{_\alpha\otimes_{\hat{\beta}}}\alpha (x)\]
\[\widehat{\Gamma}(\hat{\beta}(x))=\hat{\beta}(x)\underset{N^o}{_\alpha\otimes_{\hat{\beta}}}1\]
Following ([E2], 3.7) If $\eta_1$, $\xi_2$ are in $D(_\alpha\gH, \nu)$, let us write $(id*\omega_{\xi_2, \eta_1})(\sigma_{\nu^o}W)$ for $(\lambda_{\eta_1}^{\alpha, \hat{\beta}})^*W\rho_{\xi_2}^{\beta, \alpha}$; we have, therefore, for all $\xi_1$ and $\eta_2$ in $\gH$ :
\[(id*\omega_{\xi_2, \eta_1})(\sigma_{\nu^o}W)\xi_1|\eta_2)=(W(\xi_1\underset{\nu}{_\beta\otimes_\alpha}\xi_2)|\eta_1\underset{\nu^o}{_\alpha\otimes_{\hat{\beta}}}\eta_2)\]
Using the intertwining property of $W$ with $\alpha$, we get that it belongs to $\alpha(N)'$; we write $C_w(W)$ for the weak closure of the linear span of these operators, and we have $C_w(W)\subset \alpha(N)'$. It had been proved in ([E2], 3.10) that $C_w(W)$ is a non degenerate algebra; following ([E2] 4.1), we shall say that $W$ is weakly regular if $C_w(W)=\alpha (N)'$. If $W$ is weakly regular, then $A_w(W)=\mathcal A(W)$ and $\widehat{A_w(W)}=\widehat{\mathcal A(W)}$ ([E2], 3.12).

\subsection{Fundamental example}
\label{gd}
Let $\mathcal G$ be a measured groupo\"{\i}d, with $\mathcal G^{(0)}$ as space
of units, and $r$ and $s$ the range and source functions from $\mathcal G$ to $\mathcal G^{(0)}$, with a Haar system $(\lambda^u)_{u\in \mathcal G^{(0)}}$ and a quasi-invariant measure $\mu$ on $\mathcal G^{(0)}$. Let us write $\nu$ the associated measure on $\mathcal G$. Let us
note :
\[\mathcal G^2_{r,r}=\{(x,y)\in \mathcal G^2, r(x)=r(y)\}\]
Then, it has been shown [Val1] that the formula $W_{\mathcal G}f(x,y)=f(x,x^{-1}y)$, where $x$, $y$ are
in
$\mathcal G$, such that $r(y)=r(x)$, and $f$ belongs to $L^2(\mathcal G^{(2)})$ (with respect to an
appropriate measure, constructed from $\lambda^u$ and $\mu$), is a unitary from $L^2(\mathcal G^{(2)})$ to $L^2(\mathcal G^2_{r,r})$ (with respect also to another
appropriate measure, constructed from $\lambda^u$ and $\mu$). 
\newline
Let us define $r_{\mathcal G}$ and $s_{\mathcal G}$ from
$L^\infty (\mathcal G^{(0)})$ to $L^\infty (\mathcal G)$ (and then considered as representations on $\mathcal L(L^2(\mathcal
G))$), for any
$f$ in
$L^\infty (\mathcal G^{(0)})$, by
$r_{\mathcal G}(f)=f\circ r$ and $s_{\mathcal G}(f)=f\circ s$.
\newline
We shall identify ([Y], 3.2.2) the Hilbert space $L^2(\mathcal G^{(2)})$ with the relative Hilbert tensor product $L^2(\mathcal G, \nu)\underset{L^{\infty}(\mathcal G^{(0)}, \mu)}{_{s_{\mathcal G}}\otimes_{r_{\mathcal G}}}L^2(\mathcal G, \nu)$, and the Hilbert space $L^2(\mathcal G^2_{r,r})$ with $L^2(\mathcal G, \nu)\underset{L^{\infty}(\mathcal G^{(0)}, \mu)}{_{r_{\mathcal G}}\otimes_{r_{\mathcal G}}}L^2(\mathcal G, \nu)$. Moreover, the unitary $W_{\mathcal G}$ can be then interpreted [Val2] as a pseudo-multiplicative unitary over the basis
$L^\infty (\mathcal G^{(0)})$, with respect to the representation $r_{\mathcal G}$, and anti-representations
$s_{\mathcal G}$ and
$r_{\mathcal G}$ (as here the basis is abelian, the notions of representation and anti-representations are the same, and the commutation property is fulfilled). So, we get that $W_{\mathcal G}$ is a $(r_\mathcal G, s_\mathcal G, r_\mathcal G)$ pseudo-multiplicative unitary. 
\newline
Let us take the notations of \ref{AW}; the von Neumann algebra $\mathcal A(W_{\mathcal G})$ is equal to the von Neumann algebra $L^{\infty}(\mathcal G, \nu)$ ([Val2], 3.2.6 and 3.2.7); using ([Val2]
3.1.1), we get that the Hopf-bimodule homomorphism
$\widehat{\Gamma}$ defined on
$L^{\infty}(\mathcal G, \nu)$ by $W_{\mathcal G}$ is equal to the usual Hopf-bimodule homomorphism $\Gamma_{\mathcal G}$ studied in [Val1], and
recalled in
\ref{Hbimod}.
Moreover, the von Neumann algebra $\widehat{\mathcal A(W_{\mathcal G})}$ is equal to the von Neumann algebra $\mathcal
L(\mathcal G)$ ([Val2], 3.2.6 and 3.2.7); using ([Val2] 3.1.1), we get that the Hopf-bimodule homomorphism $\Gamma$ defined on $\mathcal L(\mathcal
G)$ by
$W_{\mathcal G}$ is the usual Hopf-bimodule homomorphism $\widehat{\Gamma_{\mathcal G}}$ studied in [Y] and [Val1]. 
\newline
Let us suppose now that the groupoid $\mathcal G$ is locally compact in the sense of [R1]; it has been proved in ([E2] 4.8) that $W_\mathcal G$ is then weakly regular (in fact was proved a much stronger condition, namely the norm regularity).

\subsection{Definitions ([L1], [L2])}
\label{LW}
Let $(N, M, \alpha, \beta, \Gamma)$ be a Hopf-bimodule, as defined in \ref{Hbimod}; a normal, semi-finite, faithful operator valued weight $T$ from $M$ to $\alpha (N)$ is said to be left-invariant if, for all $x\in \gM_T^+$, we have :
\[(id\underset{N}{_\beta*_\alpha}T)\Gamma (x)=T(x)\underset{N}{_\beta\otimes_\alpha}1\]
or, equivalently (\ref{slice}), if we choose a normal, semi-finite, faithful weight $\nu$ on $N$, and write $\Phi=\nu\circ\alpha^{-1}\circ T$, which is a normal, semi-finite, faithful weight on $M$ :
\[(id\underset{N}{_\beta*_\alpha}\Phi)\Gamma (x)=T(x)\]
A normal, semi-finite, faithful operator-valued weight $T'$ from $M$ to $\beta (N)$ will be said to be right-invariant if it is left-invariant with respect to the symmetrized Hopf-bimodule, i.e., if, for all $x\in\gM_{T'}^+$, we have :
\[(T'\underset{N}{_\beta*_\alpha}id)\Gamma (x)=1\underset{N}{_\beta\otimes_\alpha}T'(x)\]
or, equivalently, if we write $\Psi=\nu\circ\beta^{-1}\circ T'$ : 
\[(\Psi\underset{N}{_\beta*_\alpha}id)\Gamma (x)=T'(x)\]
In the case of a Hopf-bimodule, with a left-invariant normal, semi-finite, faithful operator valued weight $T$ from $M$ to $\alpha (N)$, Lesieur had constructed an isometry $U$ in the following way :
let us choose a normal, semi-finite, faithful weight $\nu$ on $N$, and let us write $\Phi=\nu\circ\alpha^{-1}\circ T$, which is a normal, semi-finite, faithful weight on $M$; let us write $H_\Phi$, $J_\Phi$, $\Delta_\Phi$ for the canonical objects of the Tomita-Takesaki theory associated to the weight $\Phi$, and let us define, for $x$ in $N$, $\hat{\beta}(x)=J_\Phi\alpha(x^*)J_\Phi$. Let $\gH$ be a Hilbert space on which $M$ is acting; then ([L2], theorem 3.14), there exists an unique isometry $U_\gH$ from $\gH\underset{\nu^o}{_\alpha\otimes_{\hat{\beta}}}H_\Phi$ to $\gH\underset{\nu}{_\beta\otimes_\alpha}H_\Phi$, such that, for any $(\beta, \nu^o)$-orthogonal 
basis $(\xi_i)_{i\in I}$ of  $(H_\Phi)_\beta$, for any $a$ in $\gN_T\cap\gN_\Phi$ and for any $v$ in $D((H_\Phi)_\beta, \nu^o)$, we have 
\[U_\gH(v\underset{\nu^o}{_\alpha\otimes_{\hat{\beta}}}\Lambda_\Phi (a))=\sum_{i\in I} \xi_i\underset{\nu}{_\beta\otimes_\alpha}\Lambda_{\Phi}((\omega_{v, \xi_i}\underset{\nu}{_\beta*_\alpha}id)(\Gamma(a)))\]
Then, Lesieur proved ([L2], theorem 3.37) that, if there exists a right-invariant normal, semi-finite, faithful operator valued weight $T'$ from $M$ to $\beta (N)$, then the isometry $U_{H_\Phi}$ is a unitary, and that $W=U_{H_\Phi}^*$ is an $(\alpha, \hat{\beta}, \beta)$-pseudo-multiplicative unitary from $H_\Phi\underset{\nu}{_\beta\otimes_\alpha}H_\Phi$ to $H_\Phi\underset{\nu^o}{_\alpha\otimes_{\hat{\beta}}}H_\Phi$\vspace{5mm}. \newline
{\bf Proposition}
{\it Let $(N, M, \alpha, \beta, \Gamma)$ be a Hopf-bimodule, as defined in \ref{Hbimod}; let us suppose that there exist  a normal, semi-finite, faithful left-invariant operator valued weight $T$ from $M$ to $\alpha (N)$ and a right-invariant normal, semi-finite, faithful operator valued weight $T'$ from $M$ to $\beta (N)$; let us write $\Phi=\nu\circ\alpha^{-1}\circ T$, and let us define, for $n$ in $N$ :
 \[\hat{\beta}(n)=J_\Phi\alpha(n^*)J_\Phi\]
Then the $(\alpha, \hat{\beta}, \beta)$-pseudo-multiplicative unitary from $H_\Phi\underset{\nu}{_\beta\otimes_\alpha}H_\Phi$ to $H_\Phi\underset{\nu^o}{_\alpha\otimes_{\hat{\beta}}}H_\Phi$ verifies, for any $x$, $y_1$, $y_2$ in $\gN_T\cap\gN_\Phi$ :}
\[(i*\omega_{J_\Phi\Lambda_\Phi (y_1^*y_2), \Lambda_\Phi (x)})(W)=
(id\underset{N}{_\beta*_\alpha}\omega_{J_\Phi\Lambda_\Phi(y_2), J_\Phi\Lambda_\Phi(y_1)})\Gamma (x^*)\]
\begin{proof}
This is just ([L2], 3.19)\vspace{5mm}. \end{proof}

{\bf Remark} Clearly, the pseudo-multplicative unitary $W$ does not depend upon the choice of the right-invariant operator-valued weight $T'$. 

\section{Coinverse and scaling group}
\label{coinverse}
In this chapter, we are dealing with a Hopf-bimodule $(N, \alpha, \beta, M, \Gamma)$, equipped with a left-invariant operator-valued weight $T_L$, and a right-invariant operator-valued weight $T_R$. If $\nu$ denotes a normal semi-finite faithful weight on the basis, let $\Phi$ (resp. $\Psi$) be the lifted normal faithful semi-finite weight on $M$ by $T_L$ (resp. $T_R$). Then, with the additional hypothesis that the two modular automorphism groups associated to the two weight $\Phi$ and $\Psi$ commute, we can construct a co-inverse, a scaling group and an antipod, using slight generalizations of the constructions made in ([L2],9) for "adapted measured quantum groupoids".

\subsection{Definition}
\label{relinv}
Let $(N, \alpha, \beta, M, \Gamma)$ be a Hopf-bimodule, equipped with a left-invariant operator-valued weight $T_L$, and a right-invariant valued weight $T_R$; let $\nu$ be a normal semi-finite faithful weight on $N$; we shall denote $\Phi=\nu\circ\alpha^{-1}\circ T_L$ and $\Psi=\nu\circ\beta^{-1}\circ T_R$ the two lifted normal semi-finite weights on $M$. We shall say that the weight $\nu$ is relatively invariant with respect to $T_L$ and $T_R$ if the two modular automorphisms groups $\sigma^\Phi$ and $\sigma^\Psi$ commute. 

\subsection{Lemma}
\label{lemE}
{\it Let $(N, \alpha, \beta, M, \Gamma)$ be a Hopf-bimodule, equipped with a left-invariant operator-valued weight $T_L$, and a right-invariant valued weight $T_R$; let $\nu$ be a normal semi-finite faithful weight on $N$, relatively invariant with respect to $T_L$ and $T_R$(\ref{relinv}); we shall denote $\Phi=\nu\circ\alpha^{-1}\circ T_L$ and $\Psi=\nu\circ\beta^{-1}\circ T_R$ the two lifted normal semi-finite weights on $M$. Let us suppose that the two modular automorphisms groups $\sigma^\Phi$ and $\sigma^\Psi$ commute, and let us denote $\delta$ the modulus of $\Psi$ with respect to $\Phi$ and $\lambda$ the scaling operator (\ref{Vaes}). We shall use the notations of \ref{propbasic}. Then :
\newline
(i) let $x\in \mathcal T_{\Psi, T_R}$ and $n\in\mathbb{N}$ and $y=e_nx$, with the notations of \ref{Vaes}; then $y$ belongs to $\gN_\Psi\cap\gN_{T_R}$, is analytical with respect to $\Psi$, and the operator $\sigma_{-i/2}^\Psi(y^*)\delta^{1/2}$ is bounded, and its closure $\overline{\sigma_{-i/2}^\Psi(y^*)\delta^{1/2}}$ belongs to $\gN_{\Phi}$; moreover, with the identifications made in \ref{Vaes}, we have :
\[\Lambda_\Phi(\overline{\sigma_{-i/2}^\Psi(y^*)\delta^{1/2}})=J_\Psi\Lambda_\Psi(y)\]
(ii) let $E$ be the linear space generated by all such elements of the form $\overline{\sigma_{-i/2}^\Psi(y^*)\delta^{1/2}}$, for all $x\in \mathcal T_{\Psi, T_R}$ and $n\in\mathbb{N}$; then $E$ is a weakly dense subspace of $\gN_\Phi$, and, for all $z\in E$, $\Lambda_\Phi(z)\in D((H_\Phi)_\beta, \nu^o)$;
\newline
(iii) the linear set of all products $<\Lambda_\Phi(z), \Lambda_\Phi(z')>_{\beta, \nu^o}$ (for $z$, $z'$ in $E$) is a dense subspace of $N$. }

\begin{proof}
As $e_n$ is analytical with respect to $\Psi$, $y$ belongs to $\gN_\Psi\cap\gN_{T_R}$, is analytical with respect to $\Psi$, and $\sigma_{-i/2}^\Psi(y^*)\delta^{1/2}$
is bounded ([V], 1.2); as $\delta^{-1}$ is the modulus of $\Phi$ with respect to $\Psi$, we get that $\overline{\sigma_{-i/2}^\Psi(y^*)\delta^{1/2}}$ belongs to $\gN_{\Phi}$; we identify $\Lambda_\Phi(\overline{\sigma_{-i/2}^\Psi (y^*)\delta^{1/2}})$ with $\Lambda_\Psi(\sigma_{-i/2}^\Psi(y^*))=J_\Psi\Lambda_\Psi(y)$, which is (i). 
\newline
The subspace $E$ contains all elements of the form $\sigma_{-i/2}^\Psi (x^*)\overline{\delta^{1/2}\sigma_{-i/2}^\Psi (e_n)}$ ($x\in\mathcal T_{\Psi, T_R}$), and, by density of $\mathcal T_{\Psi, T_R}$ in $M$, we get that the closure of $E$ contains all elements of the form $a\overline{e_n\delta^{-1/2}}\overline{\delta^{1/2}\sigma_{-i/2}^\Psi (e_n)}=ae_n\sigma_{-i/2}^\Psi(e_n)$, for all $a\in M$; now, as $e_n\sigma_{-i/2}^\Psi(e_n)$ is converging to $1$, we finally get that $E$ is dense in $M$; as $\Lambda_\Phi(E)\subset J_\Psi\Lambda_\Psi(\gN_\psi\cap\gN_{T_R})$, we get, by \ref{basic}, that, for all $z$ in $E$, $\Lambda_\Phi(z)$ belongs to $D((H_\Phi)_\beta, \nu^o)$; more precisely, we have :
\[R^{\beta, \nu^o}(\Lambda_\Phi(\sigma_{-i/2}^\Psi (x^*)\overline{\delta^{1/2}\sigma_{-i/2}^\Psi (e_n)}))=R^{\beta, \nu^o}(J_\psi\Lambda_\psi(e_nx))=\Lambda_{T_R}(e_nx)\]
Therefore, the set of elements of the form $<\Lambda_\Phi(z), \Lambda_\Phi(z')>_{\beta, \nu^o}$ contains all elements of the form $\beta^{-1}\circ T_R(x^*e_ne_nx)$, for all $x$ in $\mathcal T_{\Psi, T_R}$ and $n\in \mathbb{N}$; as $T_R(x^*e_ne_nx)=\Lambda_{T_R}(e_nx)^*\Lambda_{T_R}(e_nx)=\Lambda_{T_R}(x)^*e_n^*e_n\Lambda_{T_R}(x)$; so, its closure contains all elements of the form $\beta^{-1}\circ T_R(x^*x)$, and, therefore, it contains $\beta^{-1}\circ T_R(\gM_{T_R}^+)$, which finishes the proof. \end{proof}

\subsection{Definition}
\label{defS}
As in ([L2], 9.2), we can define, for all $\lambda\in\mathbb{C}$, a closed operator $\Delta_\Phi^\lambda\underset{N^o}{_\alpha\otimes_{\hat{\beta}}}\Delta_\Phi^\lambda$, with natural values on elementary tensor products; it is possible also to define a unitary antilinear operator $J_\Phi\underset{N^o}{_\alpha\otimes_{\hat{\beta}}}J_\Phi$ from $H_\Phi\underset{N^o}{_\alpha\otimes_{\hat{\beta}}}H_\Phi$ onto $H_\Phi\underset{N}{_{\hat{\beta}}\otimes_\alpha}H_\Phi$ (whose inverse will be $J_\Phi \underset{N}{_{\hat{\beta}}\otimes_\alpha}J_\Phi$); by composition, we define then a closed antilinear operator $S_\Phi\underset{N^o}{_\alpha\otimes_{\hat{\beta}}}S_\Phi$, with natural values on elementary tensor products, whose adjoint will be $F_\Phi \underset{N}{_{\hat{\beta}}\otimes_\alpha}F_\Phi$. 

\subsection{Proposition}
\label{propS}
{\it For all $a$, $c$ in $(\gN_\Phi\cap\gN_{T_L})^*(\gN_\Psi\cap\gN_{T_R})$, $b$, $d$ in $\mathcal T_{\Psi, T_R}$ and $g$, $h$ in $E$, the following vector :
\[U^*_{H_\Phi}\Gamma(g^*)[\Lambda_\Phi(h)\underset{\nu}{_\beta\otimes_\alpha}(\lambda^{\beta, \alpha}_{\Lambda_\Psi(\sigma_{-i}^\Psi(b^*))})^*U_{H_\Psi}(\Lambda_\Psi(a)\underset{\nu^o}{_\alpha\otimes_{\hat{\beta}}}\Lambda_\Phi((cd)^*))]\]
belongs to $D(S_\Phi\underset{\nu^*}{_\alpha\otimes_{\hat{\beta}}}S_\Phi)$, and the value of $\sigma_\nu(S_\Phi\underset{\nu^*}{_\alpha\otimes_{\hat{\beta}}}S_\Phi)$ on this vector is equal to :}
\[U^*_{H_\Phi}\Gamma(h^*)[\Lambda_\Phi(g)\underset{\nu}{_\beta\otimes_\alpha}(\lambda^{\beta, \alpha}_{\Lambda_\Psi(\sigma_{-i}^\Psi(d^*))})^*U_{H_\Psi}(\Lambda_\Psi(c)\underset{\nu^o}{_\alpha\otimes_{\hat{\beta}}}\Lambda_\Phi((ab)^*))]\]

\begin{proof}
The proof is identical to ([L2],9.9), thanks to \ref{lemE}(ii). \end{proof}

\subsection{Proposition}
\label{propG}
{\it There exists a closed densely defined anti-linear operator $G$ on $H_\Phi$ such that the linear span of :
\[(\lambda^{\beta, \alpha}_{\Lambda_\Psi(\sigma_{-i}^\Psi(b^*))})^*U_{H_\Psi}(\Lambda_\Psi(a)\underset{\nu^o}{_\alpha\otimes_{\hat{\beta}}}\Lambda_\Phi((cd)^*))\]
with $a$, $c$ in $(\gN_\Phi\cap\gN_{T_L})^*(\gN_\Psi\cap\gN_{T_R})$, $b$, $d$ in $\mathcal T_{\Psi, T_R}$, is a core of $G$, and we have :}
\begin{multline*}
G[(\lambda^{\beta, \alpha}_{\Lambda_\Psi(\sigma_{-i}^\Psi(b^*))})^*U_{H_\Psi}(\Lambda_\Psi(a)\underset{\nu^o}{_\alpha\otimes_{\hat{\beta}}}\Lambda_\Phi((cd)^*))]=
\\
(\lambda^{\beta, \alpha}_{\Lambda_\Psi(\sigma_{-i}^\Psi(d^*))})^*U_{H_\Psi}(\Lambda_\Psi(c)\underset{\nu^o}{_\alpha\otimes_{\hat{\beta}}}\Lambda_\Phi((ab)^*))
\end{multline*}

\begin{proof}
The proof is identical to ([L2],9.10), thanks to \ref{lemE}(iii). \end{proof}

\subsection{Theorem}
\label{tau}
{\it Let $(N, \alpha, \beta, M, \Gamma)$ be a Hopf-bimodule, equipped with a left-invariant operator-valued weight $T_L$, and a right-invariant valued weight $T_R$; let $\nu$ be a normal semi-finite faithful weight on $N$, relatively invariant with respect to $T_L$ and $T_R$; we shall denote $\Phi=\nu\circ\alpha^{-1}\circ T_L$ and $\Psi=\nu\circ\beta^{-1}\circ T_R$ the two lifted normal semi-finite weights on $M$. Let $G$ be the closed densely defined antilinear operator defined in \ref{propG}, and let $G=ID^{1/2}$ its polar decomposition. Then, the operator $D$ is positive self-adjoint and non singular; there exists a one-parameter automorphism group $\tau_t$ on $M$ defined, for $x\in M$, by :
\[\tau_t(x)=D^{-it}xD^{it}\]
We have, for all $n\in N$ and $t\in\mathbb{R}$ :
\[\tau_t(\alpha(n))=\alpha(\sigma^\nu_t(n))\]
\[\tau_t(\beta(n))=\beta(\sigma^\nu_t(n))\]
which allows us to define $\tau_t\underset{N}{_\beta*_\alpha}\tau_t$, $\tau_t\underset{N}{_\beta*_\alpha}\sigma_t^\Phi$ and $\sigma_t^\Psi\underset{N}{_\beta*_\alpha}\tau_{-t}$ on $M\underset{N}{_\beta*_\alpha}M$; moreover, we have :}
\[\Gamma\circ\tau_t=(\tau_t\underset{N}{_\beta*_\alpha}\tau_t)\Gamma\]
\[\Gamma\circ\sigma^\Phi_t=(\tau_t\underset{N}{_\beta*_\alpha}\sigma_t^\Phi)\Gamma\]
\[\Gamma\circ\sigma^\Psi_t=(\sigma_t^\Psi\underset{N}{_\beta*_\alpha}\tau_{-t})\Gamma\]

\begin{proof}
The proof is identical to [L2], 9.12 to 9.28. \end{proof}

\subsection{Theorem}
\label{R}
{\it Let $(N, \alpha, \beta, M, \Gamma)$ be a Hopf-bimodule, equipped with a left-invariant operator-valued weight $T_L$, and a right-invariant valued weight $T_R$; let $\nu$ be a normal semi-finite faithful weight on $N$, relatively invariant with respect to $T_L$ and $T_R$; we shall denote $\Phi=\nu\circ\alpha^{-1}\circ T_L$ and $\Psi=\nu\circ\beta^{-1}\circ T_R$ the two lifted normal semi-finite weights on $M$. Let $G$ be the closed densely defined antilinear operator defined in \ref{propG}, and let $G=ID^{1/2}$ its polar decomposition. Then, the operator $I$ is antilinear, isometric, surjective, and we have $I=I^*=I^2$; there exists a $*$-antiautomorphism $R$ on $M$ defined, for $x\in M$, by :
\[R(x)=Ix^*I\]
such that, for all $t\in\mathbb{R}$, we get $R\circ\tau_t=\tau_t\circ R$ and $R^2=id$. 
\newline
For any $a$, $b$ in $\gN_\Psi\cap\gN_{T_R}$ we have :
\[R((\omega_{J_\Psi\Lambda_\Psi(a)}\underset{N}{_\beta*_\alpha}id)\Gamma(b^*b))=
(\omega_{J_\Psi\Lambda_\Psi(b)}\underset{N}{_\beta*_\alpha}id)\Gamma(a^*a)\]
and for any $c$, $d$ in $\gN_\Phi\cap\gN_{T_L}$, we have :
\[R((id\underset{N}{_\beta*_\alpha}\omega_{J_\Phi\Lambda_\Phi(c)})\Gamma(d^*d))=
(id\underset{N}{_\beta*_\alpha}\omega_{J_\Phi\Lambda_\Phi(d)})\Gamma(c^*c))\]
\newline
For all $n\in N$, we have $R(\alpha(n))=\beta(n)$, which allows us to define $R\underset{N}{_\beta*_\alpha}R$ from $M\underset{N}{_\beta*_\alpha}M$ onto $M\underset{N^o}{_\alpha*_\beta}M$ (whose inverse will be $R\underset{N^o}{_\alpha*_\beta}R$), and we have :}
\[\Gamma\circ R=\varsigma_{N^o}(R\underset{N}{_\beta*_\alpha}R)\Gamma\]

\begin{proof}
The proof is identical to [L2], 9.38 to 9.42. \end{proof}

\subsection{Theorem}
\label{dens}
{\it Let $(N, \alpha, \beta, M, \Gamma)$ be a Hopf-bimodule, equipped with a left-invariant operator-valued weight $T_L$, and a right-invariant valued weight $T_R$; let $\nu$ be a normal semi-finite faithful weight on $N$, relatively invariant with respect to $T_L$ and $T_R$; we shall denote $\Phi=\nu\circ\alpha^{-1}\circ T_L$; then :
\newline
(i) $M$ is the weak closure of the linear span of all elements of the form $(\omega\underset{N}{_\beta*_\alpha}id)\Gamma(x)$, for all $x\in M$ and $\omega\in M_*$ such that there exists $k>0$ such that $\omega\circ\beta\leq k\nu$. 
\newline
(ii) $M$ is the weak closure of the linear span of all elements of the form $(id\underset{N}{_\beta*_\alpha}\omega)\Gamma(x)$, for all $x\in M$ and $\omega\in M_*$ such that there exists $k>0$ such that $\omega\circ\alpha\leq k\nu$. 
\newline
(iii) $M$ is the weak closure of the linear span of all elements of the form $(id*\omega_{v, w})(W)$, where $v$ belongs to $D(_\alpha H_\Phi, \nu)$ and $w$ belongs to $D((H_\Phi)_{\hat{\beta}}, \nu^o)$. }

\begin{proof}
The proof is identical to [L2], 9.25. \end{proof}

\subsection{Definition}
\label{def}
Let $(N, \alpha, \beta, M, \Gamma)$ be a Hopf-bimodule, equipped with a left-invariant operator-valued weight $T_L$, and a right-invariant valued weight $T_R$; let $\nu$ be a normal semi-finite faithful weight on $N$, relatively invariant with respect to $T_L$ and $T_R$; we shall denote $\Phi=\nu\circ\alpha^{-1}\circ T_L$ and $\Psi=\nu\circ\beta^{-1}\circ T_R$ the two lifted normal semi-finite weights on $M$; let $\tau_t$ the one-parameter automorphism group constructed in \ref{tau} and let $R$ be the involutive $*$-antiautomorphism constructed in \ref{R}. We shall call $\tau_t$ the scaling group of $(N, \alpha, \beta, M, \Gamma, T_L, T_R, \nu)$ and $R$ the coinverse of $(N, \alpha, \beta, M, \Gamma, T_L, T_R, \nu)$. Thanks to \ref{R} and \ref{dens}, we see that, $T_L$ and $\nu$ being given, $R$ does not depend on the choice of the right-invariant operator-valued weight $T_R$, provided that there exists a right-invariant operator-valued weight $T_R$ such that $\nu$ is relatively invariant with respect to $T_L$ and $T_R$. 
\newline
Similarly, from \ref{tau}, one gets that, for all $x$ in $M$, $\omega\in M_*$ such that there exists $k>0$ with $\omega\circ\alpha\leq k\nu$, $\omega'\in M_*$ such that there exists $k>0$ with $\omega\circ\beta\leq k\nu$, one has :
\[\tau_t((id\underset{N}{_\beta*_\alpha}\omega)\Gamma(x))=(id\underset{N}{_\beta*_\alpha}\omega\circ\sigma_{-t}^\Phi)\Gamma\sigma^\Phi_t(x)\]
\[\tau_t((\omega'\underset{N}{_\beta*_\alpha}id)\Gamma(x))=(\omega'\circ\sigma_{t}^\Psi\underset{N}{_\beta*_\alpha}id)\Gamma\sigma_{-t}^\Psi(x)\]
So, $T_L$ and $\nu$ being given, $\tau_t$ does not depend on the choice of the right-invariant operator-valued weight $T_R$, provided that there exists a right-invariant operator-valued weight $T_R$ such that $\nu$ is relatively invariant with respect to $T_L$ and $T_R$. 

\subsection{Theorem}
\label{thS}
{\it Let $(N, \alpha, \beta, M, \Gamma)$ be a Hopf-bimodule, equipped with a left-invariant operator-valued weight $T_L$, and a right-invariant valued weight $T_R$; let $\nu$ be a normal semi-finite faithful weight on $N$, relatively invariant with respect to $T_L$ and $T_R$; we shall denote $\Phi=\nu\circ\alpha^{-1}\circ T_L$; then, for any $\xi$, $\eta$ in $D(_\alpha H_\Phi, \nu)\cap D((H_\Phi)_{\hat{\beta}}, \nu^o)$, $(id*\omega_{\xi, \eta})(W)$ belongs to $D(\tau_{i/2})$, and, if we define $S=R\tau_{i/2}$, we have :
\[S((id*\omega_{\xi, \eta})(W))=(id*\omega_{\eta, \xi})(W)^*\]
More generally, for any $x$ in $D(S)=D(\tau_{i/2})$, we get that $S(x)^*$ belongs to $D(S)$ and $S(S(x)^*)^*=x$; 
$S$ will be called the antipod of the measured quantum groupoid, and, therefore, the co-inverse and the scaling group, given by polar decomposition of the antipod, rely only upon the pseudo-multiplicative $W$. }

\begin{proof}
It is proved similarly to [L2] 9.35 and 9.36. \end{proof}

\subsection{Proposition}
\label{RTR}
{\it Let $(N, \alpha, \beta, M, \Gamma)$ be a Hopf-bimodule, equipped with a left-invariant operator-valued weight $T_L$, and a right-invariant valued weight $T_R$; let $\nu$ be a normal semi-finite faithful weight on $N$, relatively invariant with respect to $T_L$ and $T_R$; let $\tau_t$ be the scaling group of $(N, \alpha, \beta, M, \Gamma, T_L, T_R, \nu)$ and $R$ the coinverse of $(N, \alpha, \beta, M, \Gamma, T_L, T_R, \nu)$; then :
\newline
(i) the operator-valued weight $RT_RR$ is left-invariant, the operator valued-weight $RT_LR$ is right-invariant, and $\nu$ is relatively invariant with respect to $RT_RR$ and $RT_LR$. 
\newline
(ii) $\tau_t$ is the scaling group of $(N, \alpha, \beta, M, \Gamma, RT_RR, RT_LR, \nu)$ }

\begin{proof}
Let $\Phi=\nu\circ\alpha^{-1}\circ T_L$ and $\Psi=\nu\circ\beta^{-1}\circ T_R$ the two lifted normal semi-finite weights on $M$ by $T_L$ and $T_R$;  the lifted weight by $RT_RR$ (resp. $RT_LR$) is then $\Psi\circ R$ (resp. $\Phi\circ R$). As $\sigma_t^{\Psi\circ R}=R\circ\sigma_{-t}^\Psi\circ R$ and $\sigma_s^{\Phi\circ R}=R\circ\sigma_{-s}^\Phi\circ R$, we get that $\sigma^{\Psi\circ R}$ and $\sigma^{\Phi\circ R}$ commute, which is (i). 
\newline
From \ref{tau} and \ref{R}, we get that :
\begin{multline*}
\Gamma\circ\sigma_t^{\Psi\circ R}=\Gamma\circ R\circ\sigma_{-t}^\Psi\circ R=
\varsigma_{N^o}(R\underset{N}{_\beta*_\alpha}R)\Gamma\circ\sigma_{-t}^\Psi\circ R\\
=\varsigma_{N^o}(R\circ\sigma_{-t}^\Psi\circ R\underset{N^o}{_\alpha*_\beta}R\circ\tau_t\circ R)\varsigma_{N}\Gamma=(\tau_t\underset{N}{_\beta*_\alpha}\sigma_t^{\Psi\circ R})\Gamma
\end{multline*}
from which we get that, for all $x\in M$ and $\omega\in M_*$ such that there exists $k>0$ such that $\omega\circ\alpha<k\nu$, we have :
\[\tau_t((id\underset{N}{_\beta*_\alpha}\omega)\Gamma(x))=(id\underset{N}{_\beta*_\alpha}\omega\circ\sigma_{-t}^{\Psi\circ R})\Gamma(\sigma_t^{\Psi\circ R}(x))\]
from which we get, by \ref{dens}, that $\tau_t$ is the scaling group associated to $RT_RR$, $RT_LR$ and $\nu$. \end{proof}

\section{Automorphism groups on the basis}
\label{AGB}
In this section, with the same hypothesis as in chapter \ref{coinverse}, we construct two one-parameter automorphism groups on the basis $N$ (\ref{gamma}), and we prove (\ref{thcentral}) that these automorphisms leave invariant the quasi-invariant weight $\nu$.  We prove also in \ref{thcentral} that the weight $\nu$ is also quasi-invariant with respect to $T_L$ and $RT_LR$.

\subsection{Lemma}
\label{lembeta}
{\it Let $(N, \alpha, \beta, M, \Gamma)$ be a Hopf-bimodule, equipped with a left-invariant operator-valued weight $T_L$, and a right-invariant valued weight $T_R$; let $\nu$ be a normal semi-finite faithful weight on $N$, relatively invariant with respect to $T_L$ and $T_R$. Let $x\in M\cap\alpha(N)'$ and $y\in M\cap\beta(N)'$. Then :
\newline
(i) $x$ belongs to $\beta(N)$ if and only if we have :
\[\Gamma (x)=1\underset{N}{_\beta\otimes_\alpha}x\]
(ii) $y$ belongs to $\alpha(N)$ if and only if we have :
\[\Gamma(y)=y\underset{N}{_\beta\otimes_\alpha}1\]
More generally, if $x_1$, $x_2$ are in $M\cap\alpha(N)'$ and such that $\Gamma (x_1)=1\underset{N}{_\beta\otimes_\alpha}x_2$, then $x_1=x_2\in\beta(N)$. }
\begin{proof}
The proof is given in [L2], 4.4. \end{proof}

\subsection{Proposition}
\label{gamma}
{\it Let $(N, \alpha, \beta, M, \Gamma)$ be a Hopf-bimodule, equipped with a left-invariant operator-valued weight $T_L$, and a right-invariant valued weight $T_R$; let $\nu$ be a normal semi-finite faithful weight on $N$, relatively invariant with respect to $T_L$ and $T_R$. Then, there exists a unique one-parameter group of automorphisms $\gamma_t^{L}$ of $N$ such that, for all $t\in\mathbb{R}$ and $n\in N$, we have :
\[\sigma_t^{T_L}(\beta(n))=\beta(\gamma_t^{L}(n))\]
\[\sigma_t^{RT_LR}(\alpha (n))=\alpha(\gamma_{-t}^L(n))\]
Moreover, the automorphism groups $\gamma^{L}$ and $\sigma^\nu$ commute, and there exists a positive self-adjoint non-singular operator $h_{L}$ $\eta$ $Z(N)\cap N^{\gamma^L}$ such that, for any $x\in N^+$ and $t\in \mathbb{R}$, we have :
\[\nu\circ\gamma^L_t(x)=\nu(h_{L}^tx)\]
Starting from the operator-valued weights $RT_RR$ and $RT_LR$, we obtain another one-parameter group of automorphisms $\gamma_t^R$ of $N$, such that we have :
\[\sigma_t^{RT_RR}(\beta(n))=\beta(\gamma_t^R(n))\]
\[\sigma_t^{T_R}(\alpha(n))=\alpha(\gamma_{-t}^R(n))\]
and a positive self-adjoint non-singular operator $h_R$ $\eta$ $Z(N)\cap N^{\gamma^R}$ such that we have :} 
\[\nu\circ\gamma^R_t(x)=\nu(h_R^tx)\]

\begin{proof}
The existence of $\gamma_t^{L}$ is given by [L2], 4.5; moreover, from the formula $\sigma^\Phi_t\circ\sigma^\Psi_s(\beta(n))=\sigma^\Psi_s\circ\sigma^\Phi_t(\beta(n))$, we obtain :
\[\beta(\gamma_t^{L}\circ\sigma_{-s}^\nu(n))=\beta(\sigma_{-s}^\nu\circ\gamma_t^{L}(n))\]
which gives the commutation of $\gamma_t^{L}$ and $\sigma_{-s}^\nu$. The existence of $h_{L}$ is then straightforward. The construction of $\gamma^R$ and $h_R$ is just the application of the preceeding results to $RT_RR$, $RT_LR$ and $\nu$. \end{proof}

\subsection{Proposition}
\label{proph}
{\it Let $(N, \alpha, \beta, M, \Gamma)$ be a Hopf-bimodule, equipped with a left-invariant operator-valued weight $T_L$, and a right-invariant valued weight $T_R$; let $\nu$ be a normal semi-finite faithful weight on $N$, relatively invariant with respect to $T_L$ and $T_R$. Let $T'_L$ (resp. $T'_R$) be another left (resp. right)-invariant operator-valued weight; we shall denote $\Phi=\nu\circ\alpha^{-1}\circ T_L$, $\Phi'=\nu\circ\alpha^{-1}\circ T'_L$, $\Psi=\nu\circ\beta^{-1}\circ T_R$ and $\Psi'=\nu\circ\beta^{-1}\circ T'_R$ the lifted normal semi-finite weights on $M$; then, we have :
\[\beta(h_{L}^{ist})=(D\Psi'\circ\sigma_t^\Phi:D\Psi'\circ\tau_t)_s\]
\[\alpha(h_R^{ist})=(D\Phi'\circ\sigma_{-t}^\Psi:D\Phi'\circ\tau_t)_s\]
where $\tau_s$ is the scaling group constructed from $T_L$, $T_R$ and $\nu$ as well from $RT_RR$, $RT_LR$ and $\nu$ (\ref{tau} and \ref{RTR}). }

\begin{proof}
From \ref{tau}, we get, for all $t\in\mathbb{R}$, $\Gamma\circ\sigma_t^\Phi\tau_{-t}=(id\underset{N}{_\beta*_\alpha}\sigma_t^\Phi\tau_{-t})\Gamma$, and, therefore, by the right-invariance of $T'_R$, we get, for all $x\in\gM_{T'_R}^+$, that $\tau_t\sigma_{-t}^\Phi T'_R\sigma_t^\Phi\tau_{-t}(x)= T'_R(x)$; let now $x\in\gM_{\Psi'}^+$; $T'_R(x)$ is an element of the positive extended part of $\beta(N)$ which can be written :
\[\int_0^\infty \lambda de_\lambda +(1-p)\infty\]
where $p$ is a projection in $\beta(N)$, and $e_\lambda$ is a resolution of $p$. As $x$ belongs to $\gM_{\Psi'}^+$, it is well known that $p=1$, and $T'_R(x)=\int_0^\infty \lambda de_\lambda$. There exists also a projection $q$ and a resolution of $q$ such that :
\[\tau_t\sigma_{-t}^\Phi T'_R\sigma_t^\Phi\tau_{-t}(x)=\int_0^\infty \lambda df_\lambda +(1-q)\infty\]
and, for all $\mu\in \mathbb{R}^+$, we have, because $e_\mu x e_\mu$ belongs to $\gM_{T'_R}^+$ :
\begin{eqnarray*}
e_\mu(\int_0^\infty \lambda df_\lambda)e_\mu +e_\mu(1-q)e_\mu\infty 
&=&e_\mu\tau_t\sigma_{-t}^\Phi T'_R\sigma_t^\Phi\tau_{-t}(x)e_\mu\\
&=&\tau_t\sigma_{-t}^\Phi T'_R\sigma_t^\Phi\tau_{-t}(e_\mu xe_\mu)\\
&=&T'_R(e_\mu xe_\mu)\\
&=&\int_0^\mu \lambda de_\lambda
\end{eqnarray*}
from which we infer that $(1-q)e_\mu=0$, and, therefore, that $q=1$; then, we get that $e_\mu\tau_t\sigma_{-t}^\Phi T'_R\sigma_t^\Phi\tau_{-t}(x)e_\mu$ is increasing with $\mu$ towards $T'_R(x)$. Therefore, we get 
that :
\[\tau_t\sigma_{-t}^\Phi T'_R\sigma_t^\Phi\tau_{-t}(x)\subset T'_R(x)\]
and, finally, the equality, for all $x\in\gM_{\Psi'}^+$ :
\[\tau_t\sigma_{-t}^\Phi T'_R\sigma_t^\Phi\tau_{-t}(x)=T'_R(x)\]
Moreover, as we have, for all $n\in N$ 
\[\tau_t\sigma_{-t}^\Phi(\beta(n))=\beta(\sigma_{t}^\nu\gamma_{-t}^{L}(n))\]
we get, using \ref{gamma}, that, for all $x\in \gM_{\Psi'}^+$ :
\[\Psi'(\beta(h_{L}^{-t/2})\sigma_t^\Phi\tau_{-t}(x)\beta(h_L^{-t/2}))=\Psi'(x)\]
and, therefore, that, for all $x\in M^+$ :
 \[\Psi'(\beta(h_{L}^{-t/2})\sigma_t^\Phi\tau_{-t}(x)\beta(h_L^{-t/2}))\leq\Psi'(x)\]
 A similar calculation (with $\tau_t\sigma_{-t}^\Phi$ instead of $\sigma_t^\Phi\tau_{-t}$) leads to :
 \[\Psi'(\beta(h_L^{t/2})\tau_t\sigma_{-t}^\Phi(x)\beta(h_L^{t/2}))\leq\Psi'(x)\]
 which leads to the equality, from which we get the first result. 
\newline
Applying this result to $RT_RR$, $RT_LR$ and $\nu$, we get, using again \ref{RTR} :
\begin{eqnarray*}
\beta(h_R^{ist})&=&(D\Phi'\circ R\circ\sigma_t^{\Psi\circ R}:D\Phi'\circ R\circ\tau_t)_s\\
&=&(D\Phi'\circ\sigma_{-t}^\Psi\circ R: D\Phi'\circ\tau_t\circ R)_s\\
&=&R[((D\Phi'\circ\sigma_{-t}^\Psi : D\Phi'\circ\tau_t)_{-s})^*]
\end{eqnarray*}
which leads to the result. \end{proof}

\subsection{Corollary}
\label{corh}
{\it Let $(N, \alpha, \beta, M, \Gamma)$ be a Hopf-bimodule, equipped with a left-invariant operator-valued weight $T_L$, and a right-invariant valued weight $T_R$; let $\nu$ be a normal semi-finite faithful weight on $N$, relatively invariant with respect to $T_L$ and $T_R$. We shall denote $\Phi=\nu\circ\alpha^{-1}\circ T_L$ and $\Psi=\nu\circ\beta^{-1}\circ T_R$ the two lifted normal semi-finite weights on $M$, $R$ the coinverse and $\tau_t$ the scaling group constructed in \ref{R} and \ref{tau}; we shall denote $\lambda$ the scaling operator of $\Psi$ with respect to $\Phi$ (\ref{Vaes}), $h_L$ and $h_R$ the operators constructed in \ref{gamma}. Then, for all $s$, $t$ in $\mathbb{R}$ :
\newline
(i) $(D\Psi:D\Psi\circ\tau_t)_s=\lambda^{ist}\beta(h_{L}^{ist})$
\newline
(ii) $(D\Phi:D\Phi\circ\tau_t)_s=\lambda^{ist}\alpha(h_R^{ist})$
\newline
(iii) $(D\Phi:D\Phi\circ\sigma_{-t}^{\Phi\circ R})_s=\lambda^{ist}\alpha(h_R^{ist})\alpha(h_L^{-ist})$
\newline
(iv) $(D\Psi:D\Psi\circ\sigma_t^{\Psi\circ R})_s=\lambda^{ist}\beta(h_L^{ist})\beta(h_R^{-ist})$.}
\begin{proof}
Applying \ref{proph} with $T'_R=T_R$, as $(D\Psi\circ\sigma^\Phi_t:D\Psi)_s=\lambda^{-ist}$ (\ref{Vaes}), we obtain (i). Applying \ref{proph} with $T'_L=T_L$, as $(D\Phi:D\Phi\circ\sigma_{-t}^\Psi)_s=\lambda^{ist}$, we obtain (ii). Applying \ref{proph} with $T'_R=RT_LR$, we obtain :
\begin{eqnarray*}
\beta(h_L^{ist})&=&(D\Phi\circ R\circ\sigma_t^\Phi : D\Phi\circ R\circ \tau_t)_s\\
&=&(D\Phi\circ\sigma_{-t}^{\Phi\circ R}\circ R:D\Phi\circ\tau_t\circ R)_s\\
&=&R((D\Phi\circ\sigma_{-t}^{\Phi\circ R}:D\Phi\circ\tau_t)_{-s}^*)
\end{eqnarray*}
and, therefore $\alpha(h_L^{ist})=(D\Phi\circ\sigma_{-t}^{\Phi\circ R}:D\Phi\circ\tau_t)_{-s}^*$ from which one gets :
\[\alpha(h_L^{ist})=(D\Phi\circ\sigma_{-t}^{\Phi\circ R}:D\Phi\circ\tau_t)_s\]
Using (ii), we get :
\[(D\Phi:D\Phi\circ\sigma_{-t}^{\Phi\circ R})_s=\lambda^{ist}\alpha(h_R^{ist})\alpha(h_L^{-ist})\]
which is (iii). And applying \ref{proph} with $T'_L=RT_RR$, we obtain (iv). 
\end{proof}

\subsection{Lemma}
\label{lem1tau}
{\it Let $M$ be a von Neumann algebra, $\Phi$ a normal semi-finite faithful weight on $M$, $\theta_t$ a one parameter group of automorphisms of $M$. Let us suppose that there exists a positive non singular operator $\mu$ affiliated to $M^\Phi$ such that, for all $s$, $t$ in $\mathbb{R}$, we have 
\[(D\Phi\circ\theta_t:D\Phi)_s=\mu^{ist}\]
We have then, for all $t\in\mathbb{R}$, $\theta_t(\mu)=\mu$. Let us write $\mu=\int_0^\infty \lambda de_\lambda$ the spectral decomposition of $\mu$, and let us define $f_n=\int_{1/n}^n de_\lambda$. We have then, for all $a$ in $\gN_\Phi$, $t$ in $\mathbb{R}$, $n$ in $\mathbb{N}$ :}
\[\omega_{J_\Phi\Lambda_\Phi(af_n)}\circ\theta_t=\omega_{J_\Phi\Lambda_\Phi(\theta_{-t}(a)f_n\mu^{t/2})}\]

\begin{proof}
Let us remark first that $\theta_t(\mu)=\mu$, and, therefore, $\theta_t(f_n)=f_n$. On the other hand, for any $a$ in $M$, we have :
\[\theta_{-t}\sigma_s^\Phi\theta_t(x)=\sigma_s^{\Phi\circ\theta_t}(x)=\mu^{ist}\sigma_s^\Phi(x)\mu^{-ist}\]
and then :
\[\theta_{-t}\sigma_s^\Phi(x)=\mu^{ist}\sigma_s^\Phi\theta_{-t}(x)\mu^{-ist}\]
If now $x$ is analytic with respect to $\Phi$, we get that $\theta_{-t}(f_nxf_m)$ is analytic with respect to $\Phi$ and that :
\[f_n\theta_{-t}\sigma_{i/2}^\Phi(x)f_m=\mu^{-t/2}f_n\sigma_{i/2}^\Phi(\theta_{-t}(x))f_m\mu^{t/2}\]
Let us take now $a$ in $\gN_\Phi$, analytic with respect to $\Phi$; we have, for any $y$ in $M$ :
\begin{eqnarray*}
\omega_{J_\Phi\Lambda_\Phi(f_naf_m)}\circ\theta_t(y)
&=&(\theta_t(y)J_\Phi\Lambda_\Phi(f_naf_m)|J_\Phi\Lambda_\Phi(f_naf_m))\\
&=&(\theta_t(y)\Lambda_\Phi(f_m\sigma_{-i/2}^\Phi(a^*)f_n)|\Lambda_\Phi(f_m\sigma_{-i/2}^\Phi(a^*)f_n))\\
&=&\Phi(f_n\sigma_{i/2}^\Phi(a)f_m\theta_t(y)f_m\sigma_{-i/2}^\Phi(a^*)f_n)
\end{eqnarray*}
which, using the preceeding remarks, is equal to :
\[\Phi\circ\theta_t(\mu^{-t/2}f_n\sigma_{i/2}^\Phi(\theta_{-t}(a))f_m\mu^{t/2}y\mu^{t/2}f_m\sigma_{-i/2}^\Phi(\theta_{-t}(a^*))f_n\mu^{-t/2})\]
and, making now $f_n$ increasing to $1$, we get that $\omega_{J_\Phi\Lambda_\Phi(af_m)}\circ\theta_t(y)$ is equal to :
\begin{multline*}
\Phi(\sigma_{i/2}^\Phi(\theta_{-t}(a))f_m\mu^{t/2}y\mu^{t/2}f_m\sigma_{-i/2}^\Phi(\theta_{-t}(a^*)))\\
=(y\Lambda_\Phi(f_m\mu^{t/2}\sigma_{-i/2}^\Phi(\theta_{-t}(a^*)))|\Lambda_\Phi(f_m\mu^{t/2}\sigma_{-i/2}^\Phi(\theta_{-t}(a^*)))\\
=(yJ_\Phi\Lambda_\Phi(\theta_{-t}(a)f_m\mu^{t/2})|J_\Phi\Lambda_\Phi(\theta_{-t}(a)f_m\mu^{t/2}))
\end{multline*}
from which we get the result. \end{proof}

\subsection{Lemma}
\label{lem2tau}
{\it Let $(N, \alpha, \beta, M, \Gamma)$ be a Hopf-bimodule, equipped with a left-invariant operator-valued weight $T_L$, and a right-invariant valued weight $T_R$; let $\nu$ be a normal semi-finite faithful weight on $N$, relatively invariant with respect to $T_L$ and $T_R$. We shall denote $\Phi=\nu\circ\alpha^{-1}\circ T_L$ and $\Psi=\nu\circ\beta^{-1}\circ T_R$ the two lifted normal semi-finite weights on $M$, $R$ the coinverse and $\tau_t$ the scaling group constructed in \ref{R} and \ref{tau}. Then, we have :
\newline
(i) there exists a positive non singular operator $\mu_1$ affiliated to $M^\Phi$ and invariant under $\tau_t$, such that $(D\Phi\circ\tau_t:D\Phi)_s=\mu_1^{ist}$; let us write $\mu_1=\int_0^\infty \lambda de_\lambda$ and $f_n=\int_{1/n}^nde_\lambda$; we have then, for all $a$ in $\gN_\Phi$, $t$ in $\mathbb{R}$, $n$ in $\mathbb{N}$ and $x$ in $M^+$ :
\[\omega_{J_\Phi\Lambda_\Phi(\tau_t(a)f_n)}=\omega_{J_\Phi\Lambda_\Phi(af_n\mu_1^{t/2})}\circ\tau_{-t}\]
\[T\circ\tau_t(x)=\alpha\circ\sigma_t^\nu\circ\alpha^{-1}(T(\mu_1^{t/2}x\mu_1^{-t/2}))\]
(ii) there exists a positive non singular operator $\mu_2$ affiliated to $M^\Phi$ and invariant under $\sigma_t^{\Phi\circ R}$, such that $(D\Phi\circ\sigma_{-t}^{\Phi\circ R}:D\Phi)_s=\mu_2^{ist}$; let us write $\mu_2=\int_0^\infty \lambda de'_\lambda$ and $f'_n=\int_{1/n}^nde'_\lambda$; we have then, for all $b$ in $\gN_\Phi$, $t$ in $\mathbb{R}$ and $n$ in $\mathbb{N}$ :
\[\omega_{J_\Phi\Lambda_\Phi(bf'_n)}\circ\sigma_t^{\Phi\circ R}=\omega_{J_\Phi\Lambda_\Phi(\sigma_{-t}^{\Phi\circ R}(b)f'_n\mu_2^{-t/2})}\]
\[T(\sigma_{-t}^{\Phi\circ R}(\mu_1^{-t/2}x\mu_1^{t/2}))=\alpha\circ\gamma_{t}^L\circ\alpha^{-1}(T(x))\]
Moreover, we have $\mu_1^{is}=\lambda^{-is}\alpha(h_R^{-is})$, $\mu_2^{is}=\mu_1^{is}\alpha(h_L^{is})$, and $\mu_1^{is}$, $\mu_2^{is}$, $\alpha(h_L^{is})$ belong to $\alpha(N)'\cap M^\Phi$. The non-singular operators $\mu_1$, $\mu_2$ and $\alpha(h_L)$ commute two by two. }
\begin{proof} 
By \ref{corh}(ii), we get that $(D\Phi\circ\tau_t:D\Phi)_s=\lambda^{-ist}\alpha(h_R^{-ist})$, as $\lambda$ is positive non singular, affiliated to the center $Z(M)$, and $h_R$ is positive non singular affiliated to the center of $N$, we get there exists $\mu_1$ positive non singular, affiliated to $M^\Phi$ such that :
\[\mu_1^{ist}=\lambda^{-ist}\alpha(h_R^{-ist})=(D\Phi\circ\tau_t:D\Phi)_s\]
We can then apply \ref{lem1tau} to $\tau_t$ and $\tau_t(a)f_n$ (which belongs to $\gN_\Phi$) to get the first formula of (i). On the other hand, we get that $\alpha\circ\sigma_{-t}^\nu\circ\alpha^{-1}\circ T\circ \tau_t$ is a normal semi-finite operator-valued weight which verify, for all $x\in M^+$ 
\[\alpha\circ\sigma_{-t}^\nu\circ \alpha^{-1}\circ T\circ \tau_t(x)=T(\mu_1^{t/2}x\mu_1^{t/2})\]
from which we get the second formula of (i). 
\newline
By \ref{corh}(iii), we get that $(D\Phi\circ\sigma_{-t}^{\Phi\circ R}:D\Phi)_s=\lambda^{-ist}\alpha(h_R^{-ist})\alpha(h_L^{ist})$; with the same arguments, we get that there exists $\mu_2$ positive non singular, affiliated to $M^\Phi$ such that :
\[\mu_2^{ist}=\lambda^{-ist}\alpha(h_R^{-ist})\alpha(h_L^{ist})=(D\Phi\circ\sigma_{-t}^{\Phi\circ R}:D\Phi)_s\]
and we get the first formula of (ii) by applying again \ref{lem1tau} with $\sigma_{-t}^{\Phi\circ R}$. 
\newline
On the other hand, using \ref{gamma}, we get that $\alpha\circ\gamma_{-t}^L\circ\alpha^{-1}\circ T\circ\sigma_{-t}^{\Phi\circ R}$ is an operator-valued weight which verify, for all $x\in M^+$ :
\begin{eqnarray*}
\nu\circ\alpha\circ\gamma_{-t}^L\circ\alpha^{-1}\circ T\circ\sigma_{-t}^{\Phi\circ R}(x)
&=&
\nu(h_L^{-t/2}\alpha^{-1} (T\sigma_{-t}^{\Phi\circ R}(x))h_l^{-t/2})\\
&=&\Phi(\alpha(h_L^{-t/2}\sigma_{-t}^{\Phi\circ R}(x)\alpha(h_L^{-t/2}))\\
&=&\Phi\circ\sigma_{-t}^{\Phi\circ R}[\alpha(h_L^{-t/2})x\alpha(h_L^{-t/2})]\\
&=&\Phi(\mu_2^{t/2}\alpha(h_L^{-t/2})x\alpha(h_L^{-t/2})\mu_2^{t/2})
\end{eqnarray*}
from which we get, because $\mu_2^{t/2}\alpha(h_L^{-t/2})$ commutes with $\alpha(N)$ :
\[\alpha\circ\gamma_{-t}^L\circ\alpha^{-1}\circ T\circ\sigma_{-t}^{\Phi\circ R}(x)=
T(\mu_2^{t/2}\alpha(h_L^{-t/2})x\alpha(h_L^{-t/2})\mu_2^{t/2})\]
or :
\[ T(\sigma_{-t}^{\Phi\circ R}(x))=\alpha\circ\gamma_t^L\circ\alpha^{-1} (T(\mu_1^{t/2}x\mu_1^{t/2}))\]
from which we finish the proof. \end{proof}

\subsection{Proposition}
\label{thcentral}
{\it Let $(N, \alpha, \beta, M, \Gamma)$ be a Hopf-bimodule, equipped with a left-invariant operator-valued weight $T_L$, and a right-invariant valued weight $T_R$; let $\nu$ be a normal semi-finite faithful weight on $N$, relatively invariant with respect to $T_L$ and $T_R$. We shall denote $\Phi=\nu\circ\alpha^{-1}\circ T_L$ and $\Psi=\nu\circ\beta^{-1}\circ T_R$ the two lifted normal semi-finite weights on $M$, $R$ the coinverse and $\tau_t$ the scaling group constructed in \ref{R} and \ref{tau}; let $\lambda$ be the scaling operator of $\Psi$ with respect to $\Phi$ (\ref{Vaes}), $\gamma^L$ and $\gamma^R$ the two one-parameter automorphism groups of $N$ introduced in \ref{gamma}  ; then, we have :
\newline
(i) for all $t\in\mathbb{R}$ :
\[\Gamma\circ\tau_t=(\sigma_t^{\Phi}\underset{N}{_\beta*_\alpha}\sigma_{-t}^{\Phi\circ R})\Gamma=
(\sigma_t^{\Psi\circ R}\underset{N}{_\beta*_\alpha}\sigma_{-t}^{\Psi})\Gamma\]
(ii) $h_L=h_R=1$, and :
\[\nu\circ\gamma^L=\nu\circ\gamma^R=\nu\]
(iii) for all $s$, $t$ in $\mathbb{R}$ :
\[(D\Phi:D\Phi\circ\tau_t)_s=\lambda^{ist}\]
\[(D\Psi:D\Psi\circ\tau_t)_s=\lambda^{ist}\]
(iv) for all $s$, $t$ in $\mathbb{R}$ :
\[(D\Phi\circ\sigma_t^{\Phi\circ R}: D\Phi)_s=\lambda^{ist}\]
Therefore, the modular automorphism groups $\sigma^\Phi$ and $\sigma^{\Phi\circ R}$ commute, the weight $\nu$ is relatively invariant with respect to $\Phi$ and $\Phi\circ R$ and $\lambda$ is the scaling operator of $\Phi\circ R$ with respect to $\Phi$; and we have $\tau_t(\lambda)=\lambda$, $R(\lambda)=\lambda$; 
\newline
(v) there exists a non singular positive operator $q$ affiliated to $Z(N)$ such that $\lambda=\alpha(q)=\beta(q)$. }

\begin{proof}
As, for all $n\in N$, we have :
\[\sigma_{-t}^{\Phi\circ R}(\alpha(n))=R\sigma_t^\Phi R(\alpha(n))=\alpha(\gamma_t^L(n))\]
and, by definition, $\sigma_t^\Phi (\beta(n))=\beta(\gamma_t^L(n))$, using a remark made in \ref{fiber}, we may consider the automorphism $\sigma_{-t}^\Phi\underset{N}{_\beta*_\alpha}\sigma_{t}^{\Phi\circ R}$ on $M\underset{N}{_\beta*_\alpha}M$; let's take $a$ and $b$ in $\gN_\Phi\cap\gN_{T_L}$; let's write $h_L=\int_O^\infty \lambda de^L_\lambda$ and let us write $h_p=\int_{1/p}^p de^L_\lambda$; moreover, let's use the notations of \ref{lem2tau}; we have :
\[(id\underset{N}{_\beta*_\alpha}\omega_{J_\Phi\Lambda_\Phi(b\alpha(h_p)f'_m)})(\sigma_{-t}^{\Phi}\underset{N}{_\beta*_\alpha}\sigma_{t}^{\Phi\circ R})\Gamma\circ\tau_t(f_na^*af_n)\]
is equal to :
\[\sigma_{-t}^\Phi(id\underset{N}{_\beta*_\alpha}\omega_{J_\Phi\Lambda_\Phi(b\alpha(h_p)f'_m)}\circ\sigma_t^{\Phi\circ R})\Gamma\circ\tau_t(f_na^*af_n)\]
which, thanks to \ref{lem2tau}(ii), can be written, because $\alpha(h_p)$ belongs to $\alpha(N)'\cap M^\Phi$, and therefore $b\alpha(h_p)$ belongs to $\gN_\Phi$ : 
\[\sigma_{-t}^\Phi(id\underset{N}{_\beta*_\alpha}\omega_{J_\Phi\Lambda_\Phi(\sigma_{-t}^{\Phi\circ R}(b\alpha(h_p))f'_m\mu_2^{-t/2})})\Gamma\circ\tau_t(f_na^*af_n)\]
or :
\[R\sigma_{t}^{\Phi\circ R}R(id\underset{N}{_\beta*_\alpha}\omega_{J_\Phi\Lambda_\Phi(\sigma_{-t}^{\Phi\circ R}(b\alpha(h_p))f'_m\mu_2^{-t/2})})\Gamma\circ\tau_t(f_na^*af_n)\]
By \ref{lem2tau} and \ref{lemT}, we know that $af_n\mu_1^{t/2}$ belongs to $\gN_\Phi\cap\gN_{T_L}$; using now \ref{lem2tau}(i), we get that $\tau_t(af_n)=\tau_t(a)f_n$ belongs to $\gN_\Phi\cap\gN_{T_L}$.
\newline
On the other hand, by \ref{lem2tau} and \ref{lemT}, we know that $b\alpha(h_p)f'_m$ belongs to $\gN_\Phi\cap\gN_{T_L}$; using now \ref{lem2tau}(ii), we get that :
\[\sigma_{-t}^{\Phi\circ R}(b\alpha(h_p)f'_m\mu_1^{-t/2})=\sigma_{-t}^{\Phi\circ R}(b)f'_m\mu_2^{-t/2}\alpha(h_p)\alpha(h_L^{t/2})\]
belongs to $\gN_\Phi\cap\gN_{T_L}$, and so, using again \ref{lemT}, 
\[\sigma_{-t}^{\Phi\circ R}(b)f'_m\mu_2^{-t/2}\alpha(h_p)=\sigma_{-t}^{\Phi\circ R}(b)f'_m\mu_2^{-t/2}\alpha(h_p)\alpha(h_L^{t/2})\alpha(h_p)\alpha(h_L^{-t/2})\] 
 belongs also to $\gN_\Phi\cap\gN_{T_L}$; therefore, we can use \ref{R}, and we get it is equal to : 
\[R\sigma_{t}^{\Phi\circ R}(id\underset{N}{_\beta*_\alpha}\omega_{J_\Phi\Lambda_\Phi(\tau_t(a)f_n)})\Gamma(\mu_2^{-t/2}f'_m\alpha(h_p)\sigma_{-t}^{\Phi\circ R}(b^*b)\alpha(h_p)f'_m\mu_2^{-t/2})\]
which can be written, thanks to \ref{lem2tau}(i) : 
\[R\sigma_{t}^{\Phi\circ R}(id\underset{N}{_\beta*_\alpha}\omega_{J_\Phi\Lambda_\Phi(af_n\mu_1^{t/2})}\circ\tau_{-t})\Gamma(\mu_2^{-t/2}f'_m\alpha(h_p)\sigma_{-t}^{\Phi\circ R}(b^*b)\alpha(h_p)f'_m\mu_2^{-t/2})\]
or, $\alpha(h_p)$, as well as $\mu_2^{-t/2}f'_m$, being invariant under $\sigma_t^{\Phi\circ R}$ : 
\begin{multline*}
R(id\underset{N}{_\beta*_\alpha}\omega_{J_\Phi\Lambda_\Phi(af_n\mu_1^{t/2})})(\sigma_{t}^{\Phi\circ R}\underset{N}{_\beta*_\alpha}\tau_{-t})\Gamma\circ\sigma_{-t}^{\Phi\circ R}...\\
(\mu_2^{-t/2}f'_m\alpha(h_p)b^*b\alpha(h_p)f'_m\mu_2^{-t/2})
\end{multline*}
and using \ref{tau}, and again \ref{R}, we get it is equal to : 
\begin{multline*}
R[(id\underset{N}{_\beta*_\alpha}\omega_{J_\Phi\Lambda_\Phi(af_n\mu_1^{t/2})})\Gamma(\mu_2^{-t/2}f'_m\alpha(h_p)b^*b\alpha(h_p)f'_m\mu_2^{-t/2})]\\
=
(id\underset{N}{_\beta*_\alpha}\omega_{J_\Phi\Lambda_\Phi(b\alpha(h_p)f'_m\mu_2^{-t/2})})\Gamma(\mu_1^{t/2}f_na^*af_n\mu_1^{t/2})
\end{multline*}
Finally, we have proved that, for all $a$, $b$ in $\gN_\Phi\cap \gN_{T_L}$, $m,n,p$ in $\mathbb{N}$, we have :
\begin{multline*}
(id\underset{N}{_\beta*_\alpha}\omega_{J_\Phi\Lambda_\Phi(b\alpha(h_p)f'_m)})(\sigma_{-t}^{\Phi}\underset{N}{_\beta*_\alpha}\sigma_{t}^{\Phi\circ R})\Gamma\circ\tau_t(f_na^*af_n)=\\
(id\underset{N}{_\beta*_\alpha}\omega_{J_\Phi\Lambda_\Phi(b\alpha(h_p)f'_m\mu_2^{-t/2})})\Gamma(\mu_1^{t/2}f_na^*af_n\mu_1^{t/2})
\end{multline*}
But, for all $x,y\in M$, we have :
\[\omega_{J_\Phi\Lambda_\Phi(b\alpha(h_p)f'_m)}(x)=\omega_{J_\Phi\Lambda_\Phi(b)}(\alpha(h_p)f'_mxf'_m\alpha(h_p))\]
\[\omega_{J_\Phi\Lambda_\Phi(b\alpha(h_p)f'_m\mu_2^{-t/2})}(y)=\omega_{J_\Phi\Lambda_\Phi(b)}(\alpha(h_p)f'_m\mu_2^{-t/2}x\mu_2^{-t/2}f'_m\alpha(h_p))\]
and, therefore, we get that :
\[(id\underset{N}{_\beta*_\alpha}\omega_{J_\Phi\Lambda_\Phi(b)})[(1\underset{N}{_\beta\otimes_\alpha}\alpha(h_p) f'_m)(\sigma_{-t}^{\Phi}\underset{N}{_\beta*_\alpha}\sigma_{t}^{\Phi\circ R})\Gamma\circ\tau_t(f_na^*af_n)(1\underset{N}{_\beta\otimes_\alpha}f'_m\alpha(h_p))]\]
is equal to :
\[(id\underset{N}{_\beta*_\alpha}\omega_{J_\Phi\Lambda_\Phi(b)})[(1\underset{N}{_\beta\otimes_\alpha}\alpha(h_p) f'_m\mu_2^{-t/2})\Gamma(\mu_1^{t/2}f_na^*af_n\mu_1^{t/2})
(1\underset{N}{_\beta\otimes_\alpha}\mu_2^{-t/2}f'_m\alpha(h_p))]\]
and, by density, we get that :
\[(1\underset{N}{_\beta\otimes_\alpha}\alpha(h_p) f'_m)(\sigma_{-t}^{\Phi}\underset{N}{_\beta*_\alpha}\sigma_{t}^{\Phi\circ R})\Gamma\circ\tau_t(f_na^*af_n)(1\underset{N}{_\beta\otimes_\alpha}f'_m\alpha(h_p))\]
is equal to :
\[(1\underset{N}{_\beta\otimes_\alpha}\alpha(h_p) f'_m\mu_2^{-t/2})\Gamma(\mu_1^{t/2}f_na^*af_n\mu_1^{t/2})
(1\underset{N}{_\beta\otimes_\alpha}\mu_2^{-t/2}f'_m\alpha(h_p))\]
and, after making $p$ going to $\infty$, we obtain that :
\[(1\underset{N}{_\beta\otimes_\alpha} f'_m)(\sigma_{-t}^{\Phi}\underset{N}{_\beta*_\alpha}\sigma_{t}^{\Phi\circ R})\Gamma\circ\tau_t(f_na^*af_n)(1\underset{N}{_\beta\otimes_\alpha}f'_m)\]
is equal to $(*)$:
\[(1\underset{N}{_\beta\otimes_\alpha}\ f'_m\mu_2^{-t/2})\Gamma(\mu_1^{t/2}f_na^*af_n\mu_1^{t/2}))(1\underset{N}{_\beta\otimes_\alpha}\mu_2^{-t/2}f'_m)
\]
Let's now take a file $a_i$ in $\gN_\Phi\cap\gN_{T_L}$ weakly converging to $1$; we get that 
$(1\underset{N}{_\beta\otimes_\alpha} f'_m)(\sigma_{-t}^{\Phi}\underset{N}{_\beta*_\alpha}\sigma_{t}^{\Phi\circ R})\Gamma\circ\tau_t(f_n)(1\underset{N}{_\beta\otimes_\alpha} f'_m)$ is equal to :
\[(1\underset{N}{_\beta\otimes_\alpha} f'_m\mu_2^{-t/2})\Gamma(\mu_1^{t/2}f_n\mu_1^{t/2})(1\underset{N}{_\beta\otimes_\alpha}\mu_2^{-t/2}f'_m)\]
When $n$ goes to $\infty$, then $f_n$ is increasing to $1$, the first is increasing to $1\underset{N}{_\beta\otimes_\alpha}f'_m$, and the second is increasing to :
\[(1\underset{N}{_\beta\otimes_\alpha} f'_m\mu_2^{-t/2})\Gamma(\mu_1^{t})(1\underset{N}{_\beta\otimes_\alpha}\mu_2^{-t/2}f'_m)\]
which is therefore bounded.
\newline
Taking now $m$ going to $\infty$, we get that the two non-singular operators $\Gamma(\mu_1^t)$ and $1\underset{N}{_\beta\otimes_\alpha}\mu_2^t$ are equal. Using \ref{lembeta}, we get then that $\mu_1$ is equal to $\mu_2$ (and is affiliated to $\beta(N)$), from which we get, using \ref{lem2tau}, that $h_L=1$. Applying all these calculations to $(N, \alpha, \beta, M, \Gamma, RT_RR, T_R, \nu)$, we get that $h_R=1$, which is (ii). 
\newline
Let's come back to the equality $(*)$ above; we obtain that : 
\[(1\underset{N}{_\beta\otimes_\alpha} f'_m)(\sigma_{-t}^{\Phi}\underset{N}{_\beta*_\alpha}\sigma_{t}^{\Phi\circ R})\Gamma\circ\tau_t(f_na^*af_n)(1\underset{N}{_\beta\otimes_\alpha} f'_m)\]
is equal to :
\[(1\underset{N}{_\beta\otimes_\alpha}f'_m)\Gamma(f_na^*af_n)(1\underset{N}{_\beta\otimes_\alpha} f'_m)\]
So, when $n$ and $m$ go to $\infty$, we obtain :
\[(\sigma_{-t}^{\Phi}\underset{N}{_\beta*_\alpha}\sigma_{t}^{\Phi\circ R})\Gamma\circ\tau_t(a^*a)
=\Gamma(a^*a)\]
which, by density, gives the first formula of (i), the secong being given then by \ref{RTR}. 
\newline
From (ii) and \ref{corh} (i) and (ii), we get (iii). 
\newline
From (ii) and \ref{corh}(iii), we get that $(D\Phi\circ\sigma_t^{\Phi\circ R}: D\Phi)_s=\lambda^{ist}$; therefore, as $\lambda$ is affiliated to $Z(M)$, we get the commutation of the modular groups $\sigma^\Phi$ and $\sigma^{\Phi\circ R}$. Using \ref{Vaes}, we get that there exists $\lambda_R$ positive non singular affiliated to $Z(M)$ and $\delta_R$ positive non singular affiliated to $M$ such that $(D\Phi\circ R:D\Phi)_t=\lambda_R^{it^2/2}\delta_R^{it}$, and the properties of $R$ allows us to write that $R(\lambda_R)=\lambda_R$. But, on the other hand, the formula $(D\Phi\circ\sigma_t^{\Phi\circ R}: D\Phi)_s=\lambda_R^{ist}$ (\ref{Vaes}), gives that $\lambda_R=\lambda$ and, therefore, we get that $R(\lambda)=\lambda$. The formula $\tau_t(\lambda)=\lambda$ comes from (iii), which finishes the proof of (iv). 
\newline
By (i), we have $\lambda=\mu_1=\mu_2$, and, as we had proved that $\mu_1$ is affiliated to $\beta(N)$, we get that $\lambda$ is affilated to $\beta(N)$; as $R(\lambda)=\lambda$ by (iv), we get (v). \end{proof}

\section{Measured Quantum Groupoids}
\label{MQG}

In this chapter, we give a new definition (\ref{defMQG}) of a measured quantum groupoid, and, using [L2], we get some other results, namely on the modulus (\ref{th1}), the antipod (\ref{th2}), and the manageability of the pseudo-multiplicative unitary (\ref{manageable}), all results borrowed from Lesieur. 

\subsection{Definition}
\label{defMQG}
An octuplet $(N, M, \alpha, \beta, \Gamma, T_L, T_R, \nu)$ will be called a measured quantum groupoid if :
\newline
(i) $(N, M, \alpha, \beta, \Gamma)$ is a Hopf-bimodule
\newline
(ii) $T_L$ is a normal semi-finite faithful operator-valued weight from $M$ to $\alpha(N)$, which is left-invariant, i.e. such that, for any $x\in \gM_{T_L}^+$ :
\[(id\underset{N}{_\beta*_\alpha}T_L)\Gamma(x)=T_L(x)\underset{N}{_\beta\otimes_\alpha}1\]
(iii) $T_R$ is a normal semi-finite faithful operator-valued weight from $M$ to $\beta(N)$, which is right-invariant, i.e. such that, for any $x\in \gM_{T_R}^+$ :
\[(T_R\underset{N}{_\beta*_\alpha}id)\Gamma(x)=1\underset{N}{_\beta\otimes_\alpha}T_R(x)\]
(iv) $\nu$ is a normal semi-finite faithful weight on $N$, which is relatively invariant with respect to $T_L$ and $T_R$, i.e. such that the modular automorphism groups $\sigma^\Phi$ and $\sigma^{\Psi}$ commute, where $\Phi=\nu\circ\alpha^{-1}\circ T_L$ and $\Psi=\nu\circ\beta^{-1}\circ T_R$.
\newline
Let $R$ be the co-inverse constructed in \ref{R}; thanks to \ref{thcentral}, we get that $(N, M, \alpha, \beta, \Gamma, T_L, RT_LR, \nu)$ is a measured quantum groupoid (as well as $(N, M, \alpha, \beta, \Gamma, RT_RR, T_R, \nu)$). Moreover, $R$ (resp. $\tau_t$) remains the co-inverse (resp. the scaling group) of this measured quantum groupoid. 

\subsection{Remark}
\label{rem}
Let $(N, M, \alpha, \beta, \Gamma, T_L, T_R, \nu)$ be a measured quantum groupoid in the sense of \ref{defMQG}, and let us denote $R$ (resp. $\tau_t$) the co-inverse (resp. the scaling group) constructed in \ref{R} (resp. \ref{tau}). Then $(N, M, \alpha, \beta, \Gamma, T_L, R, \tau, \nu)$ is a measured quantum groupoid in the sense of [L2], 4.1. 
\newline
Conversely if $(N, M, \alpha, \beta, \Gamma, T, R, \tau, \nu)$ is a measured quantum groupoid in the sense of [L2], 4.1, then $(N, M, \alpha, \beta, \Gamma, T, RTR, \nu)$ is a measured quantum groupoid in the sense of \ref{defMQG}. 

\subsection{Theorem}
\label{th1}
{\it Let $(N, M, \alpha, \beta, \Gamma, T_L, T_R, \nu)$ be a measured quantum groupoid; let us denote $\Phi=\nu\circ\alpha^{-1}\circ T_L$, and let $R$ be the co-inverse and $\tau_t$ the scaling group constructed in \ref{R} and \ref{tau}. Let $\delta_R$ be the modulus of $\Phi\circ R$ with respect to $\Phi$. Then, we have :
\newline
(i) $R(\delta_R)=\delta_R^{-1}$, $\tau_t(\delta_R)=\delta_R$, for all $t\in\mathbb{R}$. 
\newline
(ii) we can define a one-parameter group of unitaries $\delta_R^{it}\underset{N}{_\beta\otimes_\alpha}\delta_R^{it}$ which acts naturally on elementary tensor products, which verifies, for all $t\in\mathbb{R}$ :}
\[\Gamma(\delta_R^{it})=\delta_R^{it}\underset{N}{_\beta\otimes_\alpha}\delta_R^{it}\]

\begin{proof}
Thanks to \ref{rem}, we can rely on Lesieur's work [L2]; (i) is [L2], 5.6; (ii) is [L2], 5.20. \end{proof}

\subsection{Proposition}
\label{th2}
{\it Let $(N, M, \alpha, \beta, \Gamma, T_L, T_R, \nu)$ be a measured quantum groupoid; let us denote $\Phi=\nu\circ\alpha^{-1}\circ T_L$, and let $R$ be the co-inverse and $\tau_t$ the scaling group constructed in \ref{R} and \ref{tau}. Then :
\newline
(i) the left ideal $\gN_{T_L}\cap\gN_\Phi\cap\gN_{RT_LR}\cap\gN_{\Phi\circ R}$ is dense in $M$, and the subspace $\Lambda_\Phi(\gN_{T_L}\cap\gN_\Phi\cap\gN_{RT_LR}\cap\gN_{\Phi\circ R})$ is dense in $H_\Phi$. 
\newline
(ii) there exists a dense linear subspace $E\subset\gN_\Phi$ such that $\Lambda_\Phi(E)$ is dense in $H_\Phi$ and $J_\Phi\Lambda_\Phi(E)\subset D(_\alpha H_\Phi, \nu)\cap D((H_\Phi)_\beta, \nu^o)$. 
 }

\begin{proof}
Part (i) is given by [L2] 6.5; part (ii) by [L2] 6.7. \end{proof}

\subsection{Theorem}
\label{manageable}
{\it Let $(N, M, \alpha, \beta, \Gamma, T_L, T_R, \nu)$ be a measured quantum groupoid; let us denote $\Phi=\nu\circ\alpha^{-1}\circ T_L$, and let $R$ be the co-inverse and $\tau_t$ the scaling group constructed in \ref{R} and \ref{tau}. Then :
\newline
(i) there exists a one-parameter group of unitaries $P^{it}$ such that, for all $t\in\mathbb{R}$ and $x\in\gN_\Phi$ :
\[P^{it}\Lambda_\Phi(x)=\lambda^{t/2}\Lambda_\Phi(\tau_t(x))\]
(ii) for any $y$ in $M$, we get :
\[\tau_t(y)=P^{it}yP^{-it}\]
(iii) we have :
\[W(P^{it}\underset{N}{_\beta\otimes_\alpha}P^{it})=(P^{it}\underset{N^o}{_\alpha\otimes_{\hat{\beta}}}P^{it})W\]
(iv) for all $v\in D(P^{-1/2})$, $w\in D(P^{1/2})$, $p$, $q$ in $D(_\alpha H_\Phi, \nu)\cap D((H_\Phi)_{\hat{\beta}}, \nu^o)$, we have :
\[(W^*(v\underset{\nu^o}{_\alpha\otimes_{\hat{\beta}}}q)|w\underset{\nu}{_\beta\otimes_\alpha}p)=
(W(P^{-1/2}v\underset{\nu}{_\beta\otimes_\alpha}J_\Phi p)|P^{1/2}w\underset{\nu^o}{_\alpha\otimes_{\hat{\beta}}}J_\Phi q)\]
The pseudo-multiplicative unitary will be said to be "manageable", with "managing operator" $P$. 
\newline
(v) $W$ is weakly regular in the sense of [E2], 4.1}

\begin{proof}
The proof is given in [L2], 7.3. and 7.5. \end{proof}

\subsection{Theorem}
\label{unicity}
{\it Let $(N, M, \alpha, \beta, \Gamma, T_L, T_R, \nu)$ be a measured quantum groupoid; let us denote $\Phi=\nu\circ\alpha^{-1}\circ T_L$, and let $R$ be the co-inverse and $\tau_t$ the scaling group constructed in \ref{R} and \ref{tau}. Let $T'$ be another left-invariant operator-valued weight; let us write $\Phi'= \nu\circ\alpha^{-1}\circ T'$ and let us suppose that :
\newline
(i) $(N, M, \alpha, \beta, \Gamma, T', RT'R, \nu)$ is a measured quantum groupoid;
\newline
(ii) $\tau_t$ is the scaling group of this new quantum groupoid;
\newline
(iii) for all $t\in\mathbb{R}$, the automorphism group $\gamma^{'L}$ of $N$ defined by $\sigma_t^{\Phi'}(\beta(n))=\beta(\gamma^{'L}_t(n))$ commutes with $\gamma^L$;
\newline
Then, there exists a strictly positive operator $h$ affiliated to $Z(N)$ such that $(DT':DT)_t=\beta(h^{it})$. Moreover, we have then $\gamma^{'L}=\gamma^L$. }

\begin{proof}
This is [L2] 5.21. Then, we get :
\[\beta(\gamma^{'L}_t(n))=\sigma_t^{\Phi'}(\beta(n))=\beta(h^{-it})\beta(\gamma^L_t(n))\beta(h^{it})=\beta(\gamma_t^L(n))\]
\end{proof}


\section{Bibliography}
[BS] S. Baaj and G. Skandalis : Unitaires multiplicatifs et dualit\'{e} pour les produits crois\'{e}s de $\mathbb{C}^*$-alg\`{e}bres, {\it Ann. Sci. ENS}, {\bf 26} (1993), 425-488.
\newline\indent
[BSz1] G. B\"{o}hm and K. Szlach\'{a}nyi : A Coassociative $\mathbb{C}^*$-Quantum group with Non
Integral Dimensions, {\it Lett. Math. Phys.}, {\bf 38} (1996), 437-456.
\newline\indent
[BSz2] G. B\"{o}hm and K. Szlach\'{a}nyi : Weak $\mathbb{C}^*$-Hopf Algebras : the coassociative
symmetry of non-integral dimensions, in Quantum Groups and Quantum spaces {\it Banach Center
Publications}, {\bf 40} (1997), 9-19.
\newline\indent
[C1] A. Connes: On the spatial theory of von Neumann algebras, 
{\it J. Funct. Analysis}, {\bf 35} (1980), 153-164.
\newline\indent
[C2] A. Connes: Non commutative Geometry, Academic Press, 1994
\newline\indent
[E1] M. Enock : Inclusions irr\'{e}ductibles de facteurs et unitaires multiplicatifs II, {\it J. Funct. Analysis}, {\bf 137} (1996), 466-543. 
\newline\indent
[E2] M. Enock : Quantum groupoids of compact type, {\it J. Inst. Math. Jussieu}, {\bf 4} (2005), 29-133. 
\newline\indent
[E3] M. Enock : Inclusions of von Neumann algebras and quantum groupo\"{\i}ds III,  {\it J. Funct. Analysis}, {\bf 223} (2005), 311-364. 
\newline\indent
[EN] M. Enock, R. Nest : Inclusions of factors, multiplicative unitaries and Kac algebras, {\it J. Funct. Nanlysis}, {\bf 137} (1996), 466-543. 
\newline\indent
[ES] M. Enock, J.-M. Schwartz : Kac Algebras and Duality of locally compact Groups, Springer-Verlag, Berlin, 1989. 
\newline\indent
[EV] M. Enock, J.-M. Vallin : Inclusions of von Neumann algebras and quantum groupo\"{\i}ds,
{\it J. Funct. Analalysis}, {\bf 172} (2000), 249-300.
\newline\indent
[J] V. Jones: Index for subfactors, {\it Invent. Math.},  {\bf 72} (1983), 1-25.
\newline\indent
[KV1] J. Kustermans and S. Vaes : Locally compact quantum groups, {\it Ann. Sci. ENS}, {\bf 33} (2000), 837-934.
\newline\indent
[KV2] J. Kustermans and S. Vaes, Locally compact quantum groups in the von Neumann algebraic setting, {\it Math. Scand.}, {\bf 92} (2003), 68-92. 
\newline\indent
[L1] F. Lesieur : thesis, University of Orleans, available at :
\newline
http://tel.ccsd.cnrs.fr/documents/archives0/00/00/55/05
\newline\indent
[L2] F. Lesieur : Measured Quantum groupoids, math.OA/0504104, to be published in {\it M\' emoires de la SMF} (2007).  
\newline\indent
[MN] T. Masuda and Y. Nakagami : A von Neumann Algebra framework for the duality of the quantum groups, {\it Publ. RIMS Kyoto}, {\bf 30} (1994), 799-850.
\newline\indent
[MNW] T. Masuda, Y. Nakagami and S.L. Woronowicz : A $\mathbb{C}^*$-algebraic framework for quantum groups, {\it Internat. J. math. }, {\bf 14} (2003), 903-1001. 
\newline\indent
[NV1] D. Nikshych, L. Va\u{\i}nerman : Algebraic versions of a finite dimensional quantum
groupoid, in Lecture Notes in Pure and Applied Mathematics, Marcel Dekker, 2000.
\newline\indent
[NV2] D. Nikshych, L. Va\u{\i}nerman : A characterization of depth 2 subfactors of $II_1$
factors, {\it J. Funct. Analysis}, {\bf 171} (2000), 278-307.
\newline\indent
[R1] J. Renault : A Groupoid Approach to $\mathbb{C}^*$-Algebras, {\it Lecture Notes in
Math.} {\bf 793}, Springer-Verlag
\newline\indent
[R2] J. Renault : The Fourier algebra of a measured groupoid and its multipliers, {\it J. Funct. Analysis}, {\bf 145} (1997), 455-490.
\newline\indent
[S1] J.-L. Sauvageot : Produit tensoriel de $Z$-modules et applications, in Operator
Algebras and their Connections with Topology and Ergodic Theory, Proceedings Bu\c{s}teni, Romania,
1983, {\it Lecture Notes in
Math.} {\bf 1132}, Springer-Verlag, 468-485.
\newline\indent
[S2] J.-L. Sauvageot : Sur le produit tensoriel
relatif d'espaces de Hilbert,  {\it J. Operator Theory}, {\bf 9} (1983), 237-352.
\newline\indent
[Sz] K. Szlach\'{a}nyi : Weak Hopf algebras, in Operator Algebras and Quantum Field Theory,
S. Doplicher, R. Longo, J.E. Roberts, L. Zsido editors, International Press, 1996.
\newline\indent
[T] M. Takesaki : Theory of Operator Algebras II, Springer, Berlin, 2003. 
\newline\indent
[V] S. Vaes : A Radon-Nikodym theorem for von Neumann algebras, {\it J. Operator Theory}, {\bf 46} (2001), 477-489.
\newline\indent
[Val1] J.-M. Vallin : Bimodules de Hopf
et Poids op\'eratoriels de Haar, {\it J. Operator theory}, {\bf 35} (1996), 39-65
\newline\indent
[Val2] J.-M. Vallin : Unitaire pseudo-multiplicatif associ\'e \`a un groupo\"{\i}de;
applications \`a la moyennabilit\'e,  {\it J. Operator theory}, {\bf 44} (2000), 347-368.
\newline\indent
[Val3] J.-M. Vallin : Groupo\"{\i}des quantiques finis, {\it J. Algebra}, {\bf 239} (2001), 215-261.
\newline\indent
[Val4] J.-M. Vallin : Multiplicative partial isometries and finite quantum groupoids, in Locally Compact Quantum Groups and Groupoids, IRMA Lectures in Mathematics and Theoretical Physics 2, V. Turaev, L. Vainerman editors, de Gruyter, 2002. 
\newline\indent
[W1] S.L. Woronowicz : Tannaka-Krein duality for compact matrix pseudogroups. Twisted $SU(N)$ group. {\it Invent. Math.}, {\bf 93} (1988), 35-76. 
\newline\indent
[W2] S.L. Woronowicz : Compact quantum group, in "Sym\'{e}tries quantiques" (Les Houches, 1995), North-Holland, Amsterdam (1998), 845-884.
\newline\indent
[W3] S.L. Woronowicz : From multiplicative unitaries to quantum groups, {\it Int. J. Math.}, {\bf 7}
(1996), 127-149. 
\newline\indent
[Y] T. Yamanouchi : Duality for actions and coactions of
measured Groupoids on von Neumann Algebras, {\it Memoirs of the A.M.S.}, {\bf 101} (1993), 1-109.
\newline\indent

\end{document}